\documentclass[11pt]{amsart}
\usepackage{amsmath, amsfonts, amssymb, color, enumitem}

\usepackage{graphicx}
\usepackage{geometry}
\usepackage[utf8]{inputenc}
\usepackage[T1]{fontenc}
\usepackage{xcolor}
\usepackage{hyperref}

\usepackage{ulem}

\usepackage{tikz}
\usetikzlibrary{arrows}
\usepackage{pgfplots}
\usetikzlibrary{patterns}

\usepackage{dsfont}

\definecolor{boubelcolor}{rgb}{.65,0.05,0}

\DeclareUnicodeCharacter{2260}{\colorlet{memoireboubel}{.}\color{boubelcolor}} 
\DeclareUnicodeCharacter{00B1}{\color{memoireboubel}{}} 

\DeclareUnicodeCharacter{00A3}{\colorlet{memoirejuillet}{.}\color{blue}} 
\DeclareUnicodeCharacter{00A7}{\color{memoirejuillet}{}} 

\newtheorem{theorem}{Theorem}[section]
\newtheorem{lemma}[theorem]{Lemma}
\newtheorem{proposition}[theorem]{Proposition}
\newtheorem{corollary}[theorem]{Corollary}
\newtheorem{definition}[theorem]{Definition}
\newtheorem{remark}[theorem]{Remark}

\numberwithin{equation}{section}

\newcommand{\de}{\, \mathrm{d}} 
\newcommand{\Law}{\mathrm{Law}}

\newcommand{\shadow}[2]{\mathcal{S}^{#1}(#2)}

\newcommand{\leqc}{\leq_{c}}
\newcommand{\leqp}{\leq_{+}}

\newcommand{\MO}{\mathcal{M}_1}

\newcommand{\PO}{\mathcal{P}_1}

\newcommand{\TO}{\mathcal{T}_1}
\newcommand{\1}{\mathds{1}}
\newcommand{\SEP}{\mathrm{SEP}}

\newcommand{\R}{\mathbb{R}}
\renewcommand{\P}{\mathbb{P}}

\begin{document}
	\title[Shadows and Barriers]{Shadows and Barriers}
	\author{Martin Br\"uckerhoff \address[Martin Br\"uckerhoff]{Universit\"at M\"unster, Germany} \email{martin.brueckerhoff@uni-muenster.de} \hspace*{0.5cm} Martin Huesmann \address[Martin Huesmann]{Universit\"at M\"unster, Germany} \email{martin.huesmann@uni-muenster.de} } 
	\thanks{MB and MH are funded by the Deutsche Forschungsgemeinschaft (DFG, German Research Foundation) under Germany's Excellence Strategy EXC 2044 –390685587, Mathematics Münster: Dynamics–Geometry–Structure. }
	
	\begin{abstract}
		We show an intimate connection between solutions of the Skorokhod Embedding Problem which are given as the first hitting time of a barrier and the concept of shadows in martingale optimal transport. More precisely, we show that a solution $\tau$ to the Skorokhod Embedding Problem between $\mu$ and $\nu$ is of the form $\tau = \inf \{t \geq 0 : (X_t,B_t) \in \mathcal{R}\}$ for some increasing process $(X_t)_{t \geq 0}$ and a barrier $\mathcal{R}$  if and only if there exists a  time-change $(T_l)_{l \geq 0}$ such that for all $l \geq 0$ the equation
		$$\P[B_{\tau} \in \cdot , \tau \geq T_l] = \shadow{\nu}{\P[B_{T_l} \in \cdot , \tau \geq T_l]}$$
		is satisfied, i.e.\ the distribution of $B_{\tau}$ on the event that the Brownian motion is stopped after $T_l$ is  the shadow of the distribution of $B_{T_l}$ on this event in the terminal distribution $\nu$.
		
		This equivalence allows us to construct new families of barrier solutions that naturally interpolate between two given barrier solutions. We exemplify this by an interpolation between the Root embedding and the left-monotone embedding.

		\smallskip
		
		\normalem

	\end{abstract}
	
	\date{\today}
	\maketitle

	\section{Introduction}
	
	Let $\mu$ be a probability measure on $\mathbb{R}$ and $(B_t)_{t \geq 0}$ a $\mathcal{F}$-Brownian Motion with initial distribution $\mathrm{Law}(B_0) = \mu$ defined on a filtered probability space $(\Omega, \mathcal{A},\mathbb{P}, (\mathcal{F}_t)_{t \geq 0})$. We assume that the filtration $\mathcal{F} = (\mathcal{F}_t)_{t \geq 0}$ is right-continuous and completed w.r.t.\ $\mathbb{P}$.
	
	Given another probability measure $\nu$, a finite $\mathcal{F}$-stopping time $\tau$ is said to be a solution to the Skorokhod Embedding Problem w.r.t.\ $\mu$ and $\nu$, if 
	\begin{equation}
	\tag{$\mathrm{SEP}(\mu,\nu)$}
	(B_{t \land \tau})_{t \geq 0} \text{ is uniformly integrable} \quad \text{and} \quad B_{\tau} \sim \nu.
	\end{equation}
	It is well known that there exists a solution to $\SEP(\mu,\nu)$ if and only if $\mu \leqc \nu$, i.e.\ if we have $\int_{\mathbb{R}} \varphi \de \mu \leq \int _{\mathbb{R}} \varphi \de \nu$ for all convex functions $\varphi$. In general there exist many different solutions to $\SEP(\mu,\nu)$ (cf.\ \cite{Ob04}).

	\subsection*{Main Result}
	
	In this article, we focus on the subclass of ``barrier solutions'' to the Skorokhod Embedding Problem which includes for instance the Root embedding \cite{Ro69}, the Az\'{e}ma-Yor embedding \cite{AzYo79}, the Vallois embedding \cite{Va83}, and the left-monotone embedding \cite{BeHeTo17}. These solutions can be described as the first time the process $(X_t,B_t)_{t \geq 0}$ hits a barrier in $[0, \infty) \times \mathbb{R}$ (cf.\ Definition  \ref{def:Intro})  where $X$ is monotonously increasing, non-negative and $\mathcal{F}$-adapted. 
	We show that  these embeddings are closely related to the concept of shadows introduced by Beiglb\"ock and Juillet in \cite{BeJu16}.
	
	\begin{definition} \label{def:Intro}
		\begin{itemize}
			\item [(i)]  A set $\mathcal{R} \subset [0, \infty) \times \mathbb{R}$ is called a barrier if $\mathcal{R}$ is closed and for all $(l,x) \in \mathcal{R}$ and $l \leq l'$ we have $(l',x) \in \mathcal{R}$.
			
			\item [(ii)] Let $\xi$ and $\zeta$ be finite measures on $\mathbb{R}$. We say that $\xi$ is a submeasure of $\zeta$ if $\xi[A] \leqp \zeta[A]$ for all $A \in \mathcal{B}(\mathbb{R})$, denoted by $\xi \leqp \zeta$. 
			
			\item [(iii)] Let $\eta$ and $\zeta$ be  finite measures on $\mathbb{R}$. A finite measure $\xi$ that satisfies $\eta \leqc \xi \leqp \zeta$ and $\xi \leqc \xi'$ for all $\xi'$ with $\eta \leqc \xi' \leqp \zeta$, is called the shadow of $\eta$ in $\zeta$ and is denoted by $\shadow{\zeta}{\eta}$.
		\end{itemize}
	\end{definition}
	
	We want to mention that in the literature  barriers defined as in Definition \ref{def:Intro}(i)   are  sometimes  called ``right-barriers'' in contrast to ``left-barriers''.
	The shadow $\shadow{\zeta}{\eta}$ exists whenever the set of possible candidates is not empty, i.e.\ if there exists $\xi$ such that $\eta \leqc \xi \leqp \zeta$.
	This existence result was first shown by Rost \cite{Ro71}. Later Beiglb\"ock and Juillet \cite{BeJu16} rediscovered this object in the context of martingale optimal transport and coined the name shadow. 
	
	In the following we use the notation $\Law(X;A)$ for the (sub-)probability measure which is given by the push-forward of $X$ under the restriction of $\mathbb{P}$ to the event $A$ (cf.\ Section \ref{ssec:Notation}).

	\begin{theorem} \label{thm:intro}
		Let $\mu \leqc \nu$ and $\tau$ a solution of $\mathrm{SEP}(\mu,\nu)$.
		The following are equivalent:
		\begin{itemize}
			\item [(i)] There exists a right-continuous $\mathcal{F}$-adapted stochastic process $(X_t)_{t \geq 0}$ which is non-negative, monotonously increasing and satisfies $\mathbb{P}[\exists s < t : X_s = X_{t} = l] = 0$ for all $l \geq 0$, and a closed barrier $\mathcal{R} \subset [0, \infty) \times \mathbb{R}$ such that 
			\begin{equation*}
			\tau = \inf \{ t \geq 0 : (X_t,B_t) \in \mathcal{R}\} \quad a.s.
			\end{equation*}
			
			\item [(ii)] There exists a left-continuous $\mathcal{F}$-time-change $(T_l)_{l \geq 0}$ with $T_0 = 0$, $T_\infty = \infty$ and $\mathbb{P}[\lim _{k \downarrow l} T_k = T_l] = 1$ for all $l \geq 0$ such that for all $l \geq 0$ we have
			\begin{equation} \label{eq:ShadowResid}
			\mathrm{Law}(B_{\tau}; \tau \geq T_l) = \shadow{\nu}{\mathrm{Law}(B_{T_l}; \tau \geq T_l)}.
			\end{equation}
		\end{itemize}
		Moreover, we may choose $T_l = \inf \{t \geq 0: X_t \geq l\}$ or $X_t := \sup \{l \geq 0: T_l \leq t\}$ and $\mathcal{R} := \{ (l,x) \in [0, \infty) \times \mathbb{R} : U_{\Law(B_{T_l \land \tau})} (x) = U_{\nu}(x) \}$, respectively, where $U_\cdot$ denotes the potential function of a finite measure (see Definition \ref{def:PotentialFunction}).
	\end{theorem}
	
	\begin{remark}
		We want to stress that this theorem also holds for randomized stopping times (cf.\ Theorem \ref{thm:MainEqui}). This concerns the implication $(ii) \Rightarrow (i)$, as part (i) already ensures that the randomized stopping time is induced by a (non-randomized) stopping time.
	\end{remark}
	
	To the best of our knowledge, the only known connection between shadows and the Skorokhod Embedding Problem is implicitly through the left-monotone embedding because it is uniquely characterized by the property that the induced martingale coupling between the initial and the terminal marginal distribution is precisely the left-curtain coupling (see below). Theorem \ref{thm:intro} shows that this connection is not by accident, but just a special case of an intimate connection between shadows and barrier solutions rooted in potential theory.

	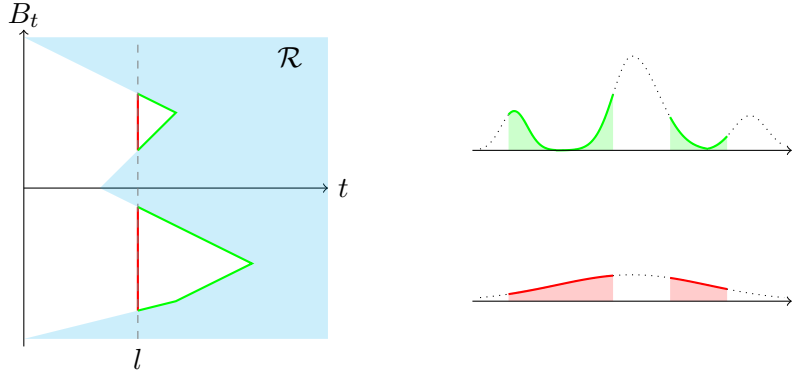
\begin{figure}  
		\begin{center}
			\begin{tikzpicture}
			\draw[->] (0,-2.1) -- (0,2.1);
			\node[above] at (0,2) {$B_t$};
			\draw[->] (0,0) -- (4,0);
			\node[right] at (4,0) {$t$};
			
			\fill[fill = cyan, opacity = 0.2] (1,0) -- (2,1) -- (0,2) -- (4,2) -- (4,0) -- (1,0);
			\fill[fill = cyan, opacity = 0.2] (0,-2) -- (2,-1.5) -- (3,-1) -- (1,0) -- (4,0) -- (4,-2)  -- (0,-2);
			\node[below] at (3.5,2) {$\mathcal{R}$};
			
			\draw[color = red,thick] (1.5,0.5) -- (1.5,1.25);
			\draw[color = green,thick] (1.5,0.5) -- (2,1) -- (1.5,1.25);
			\draw[color = red,thick] (1.5,-0.25) -- (1.5,-1.625);
			\draw[color = green,thick] (1.5,-0.25) -- (3,-1) -- (2,-1.5) -- (1.5,-1.625);
			\draw[dashed,thin,color=gray] (1.5,-2) -- (1.5,2);
			\node[below] at (1.5,-2) {$l$};

			\fill[fill = red, opacity = 0.2] (6.375,-1.5) -- (6.375,-1.4) -- (6.5,-1.37) --  (7.75,-1.15) -- (7.75,-1.5);
			\fill[fill = red, opacity = 0.2] (8.5,-1.5) -- (8.5,-1.2) -- (9,-1.3) --  (9.25,-1.35) -- (9.25,-1.5);
			\draw[dotted, domain = 6:10,smooth,variable=\t] 
			plot({\t},{0.35*exp(-(\t-8)*(\t-8)/2)-1.5});
			\draw[thick,red, domain = 6.375:7.75,smooth,variable=\t] 
			plot({\t},{0.35*exp(-(\t-8)*(\t-8)/2)-1.5});
			\draw[thick,red, domain = 8.5:9.25,smooth,variable=\t] 
			plot({\t},{0.35*exp(-(\t-8)*(\t-8)/2)-1.5});

			\draw[dotted, domain = 6:7,smooth,variable=\t] 
			plot({\t},{0.5*exp(-(\t-6.4)*(\t-6.5)/0.06)+0.5});
			\draw[dotted, domain = 7:8,smooth,variable=\t] 
			plot({\t},{1.25*exp(-(\t-8)*(\t-8)/0.12)+0.5});
			\draw[dotted, domain = 8:9,smooth,variable=\t] 
			plot({\t},{1.25*exp(-(\t-8)*(\t-8)/0.24)+0.5});
			\draw[dotted, domain = 9:10,smooth,variable=\t] 
			plot({\t},{0.45*exp(-(\t-9.6)*(\t-9.5)/0.1)+0.5});
			
			\fill[fill = green, opacity = 0.2] (6.375,0.5) -- (6.375,1) -- (6.45,1.02) -- (6.6,0.9) -- (6.75,0.6) -- (6.9,0.5);
			\fill[fill = green, opacity = 0.2] (7.2,0.5) -- (7.3,0.5)-- (7.5,0.63) -- (7.6,0.8) -- (7.75,1.28) -- (7.75,0.5);
			\fill[fill = green, opacity = 0.2] (8.5,0.5) -- (8.5,0.95) -- (8.7,0.65) -- (8.85,0.55) -- (9,0.5);
			\fill[fill = green, opacity = 0.2] (9,0.5) -- (9.25,0.7) --  (9.25,0.5);
			\draw[thick,green, domain = 6.375:7,smooth,variable=\t] 
			plot({\t},{0.5*exp(-(\t-6.4)*(\t-6.5)/0.06)+0.5});
			\draw[thick,green, domain = 7:7.75,smooth,variable=\t] 
			plot({\t},{1.25*exp(-(\t-8)*(\t-8)/0.12)+0.5});
			\draw[thick,green, domain = 8.5:9,smooth,variable=\t] 
			plot({\t},{1.25*exp(-(\t-8)*(\t-8)/0.24)+0.5});
			\draw[thick, green, domain = 9:9.25,smooth,variable=\t] 
			plot({\t},{0.45*exp(-(\t-9.6)*(\t-9.5)/0.1)+0.5});
			
			\draw[->] (5.9,0.5) -- (10.1,0.5);
			\draw[->] (5.9,-1.5) -- (10.1,-1.5);
			\end{tikzpicture}
			\caption{This is a sketch of the support and densities of the measures appearing in \eqref{eq:RootExpl}.	The supports of $\Law((l,B_l); \tau \geq l)$ (red) and $\Law((\tau,B_\tau); \tau \geq l)$ (green) are shown on the left and the densities of $\Law(B_l; \tau \geq l)$ (red) and $\Law(B_\tau; \tau \geq l)$ (green) on the right.} 
			\label{fig:RootExp}
		\end{center}
	\end{figure}
	
	\begin{figure} 
		\begin{center}		

			\begin{tikzpicture}
			\draw[->] (0,-2.1) -- (0,2.1);
			\node[above] at (0,2) {$B_t$};
			\draw[->] (0,0) -- (4,0);
			\node[right] at (4,0) {$e^{-B_0}$};
			
			\draw[dotted, domain = 0.14:4,smooth,variable=\t] 
			plot({\t},{-ln(\t)});

			\fill[fill = cyan, thick, opacity = 0.2] (0,2) --(0,1.5) -- (0.368,1.5) -- (0.5,1.234) -- (0.6,1.091) -- (0.75,0.932) -- (1,0.75) -- (1.25,0.626) -- (1.649,0.5) -- (1.649,-1) --  (0.368,-1.5) -- (0,-1.5) -- (0,-2) -- (4,-2) -- (4,2);
			\draw[cyan,opacity = 0.2, domain = 0.368:1.649,smooth,variable=\t] 
			plot({\t},{1/6*ln(\t)*ln(\t) - 7/12 *ln(\t) +3/4});
			\draw[cyan, opacity = 0.2, thick, domain = 0.368:1.649,smooth,variable=\t] 
			plot({\t},{1/6*ln(\t)*ln(\t)+5/12*ln(\t)-5/4});

			\node[below] at (3.5,2) {$\mathcal{R}$};

			\draw[red, thick,domain = 1:1.649,smooth,variable=\t] 
			plot({\t},{-ln(\t)});
			\draw[blue, thick,domain = 1.649:2.718,smooth,variable=\t] 
			plot({\t},{-ln(\t)});
			\draw[green, thick, domain = 1:1.649,smooth,variable=\t] 
			plot({\t},{1/6*ln(\t)*ln(\t)+5/12*ln(\t)-5/4});
			\draw[green, thick, domain = 1:1.649,smooth,variable=\t] 
			plot({\t},{1/6*ln(\t)*ln(\t) - 7/12 *ln(\t) +3/4});
			\draw[dashed,thin,color=gray] (1,-2) -- (1,2);
			\node[below] at (1,-2) {$l$};

			\draw[->] (5.9,0.5) -- (10.1,0.5);
			\draw[dotted] (7,-1.5) -- (7,-1) -- (9,-1) -- (9,-1.5);
			\fill[fill = blue, opacity = 0.2] (7,-1.5) -- (7,-1) -- (7.5,-1) -- (7.5,-1.5);
			\draw[color = blue,thick] (7,-1) --  (7.5,-1);
			\fill[fill = red, opacity = 0.2] (7.5,-1.5) -- (7.5,-1) -- (8,-1) -- (8,-1.5);
			\draw[color = red,thick] (7.5,-1) --  (8,-1);
			
			\draw[->] (5.9,-1.5) -- (10.1,-1.5);
			\draw[dotted] (6.5,0.5) -- (6.5,1) -- (7.5,1) -- (7.5,0.5);
			\draw[dotted] (8.5,0.5) -- (8.5,1) -- (9.5,1) -- (9.5,0.5);
			
			\fill[fill = green, opacity = 0.2] (6.75,0.5) -- (6.75,1) -- (7,1) -- (7,0.5);
			\draw[thick, color = green] (6.8,1) -- (7,1);
			\fill[fill = blue, opacity = 0.2] (7,0.5) -- (7,1) -- (7.5,1) -- (7.5,0.5);
			\draw[color = blue,thick] (7,1) --  (7.5,1);
			\fill[fill = green, opacity = 0.2] (8.5,0.5) -- (8.5,1) -- (8.75,1) -- (8.75,0.5);
			\draw[thick, color = green] (8.5,1) -- (8.8,1);
			
			\end{tikzpicture}
			
			\caption{This is a sketch of the support and densities of the measure appearing in \eqref{eq:LMExpl}.	The supports of $\Law((l,B_0); B_0 \leq - \ln(l), \tau > 0)$ (red),  $\Law((l,B_0); B_0 \leq - \ln(l), \tau = 0)$ (blue) and $\Law((\tau,B_{\tau}); B_0 \leq - \ln(l), \tau > 0)$ (green) are shown on the left and the densities of $\Law(B_0; B_0 \leq - \ln(l), \tau > 0)$ (red),  $\Law((B_0; B_0 \leq - \ln(l), \tau = 0)$ (blue) and $\Law(B_{\tau}; B_0 \leq - \ln(l), \tau > 0)$ on the right.}
			\label{fig:LmExp}
		\end{center}
	\end{figure}
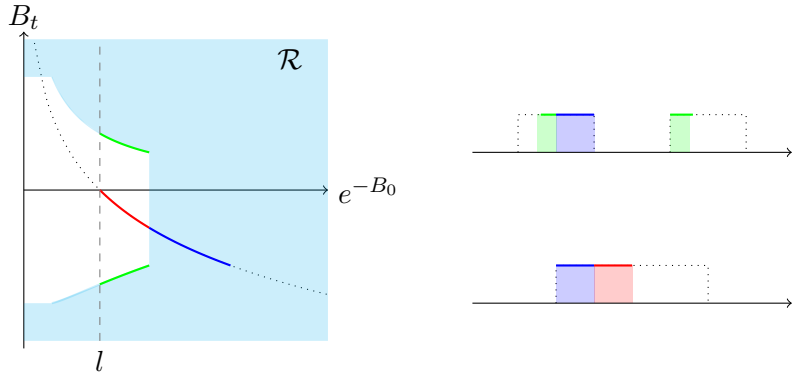

	Since the time change in Theorem \ref{thm:intro} is given by $T_l = \inf \{t \geq 0 : X_t \geq l\}$ it is straightforward to compute the time changes for well known examples.
	In the case of the Root-embedding,  we have $X^{r}_t := t$, $T^r_l = l$ and Property \eqref{eq:ShadowResid} turns into 
	\begin{equation}\label{eq:RootExpl}
	\Law(B_\tau; \tau \geq l) = \shadow{\nu}{\Law(B_l; \tau \geq l)}
	\end{equation}
	for all $l \geq 0$. 
	The measure $\Law(B_\tau; \tau \geq l)$ is the projection of $\Law((\tau,B_{\tau}); \tau \geq l)$ onto the second (the spatial) component. In the SEP context the joint law of $(\tau,B_\tau)$ describes when and where the Brownian motion is stopped. Since $\tau$ is a barrier stopping time the support of $\Law((\tau,B_{\tau}); \tau \geq l)$ is on the boundary  of the barrier intersected with $[l,\infty)\times \R$. This is depicted on the left hand side of Figure \ref{fig:RootExp}. By \eqref{eq:RootExpl} we can characterize this measure using information from time $l$ only. For each $l \geq 0$, it is given 
	as the shadow of $\Law(B_l; \tau \geq l)$ in the prescribed terminal distribution $\nu$.
	
	We have a similar situation in the case of the left-monotone embedding. We have $X^{lm}_t = \exp(-B_0)$, 
	\begin{equation*}
	T^{lm}_l = \begin{cases}
	0 & \exp(-B_0) \geq l \\ 
	+ \infty &\exp(-B_0) < l
	\end{cases}
	\end{equation*}
	and Property \eqref{eq:ShadowResid} becomes 
	\begin{equation} \label{eq:LMExpl}
	\Law(B_{\tau}; B_0 \leq -\ln(l)) = \shadow{\nu}{\mathrm{Law}(B_0; B_0 \leq - \ln(l) )}
	\end{equation}
	for all $l \geq 0$. Again the measure $\Law(B_\tau; \tau \geq l)$ is the projection of $\Law((\tau,B_{\tau}); \tau \geq l)$ onto the second component and in the SEP context the latter measure is supported on the boundary of the barrier after time $l$ (left side of Figure \ref{fig:LmExp}). Recall that in the left-monotone phase space, the Brownian motion is only moving vertically. This time the characterization of $\Law(B_\tau; \tau \geq l)$ via the shadow of $\mathrm{Law}(B_0; B_0 \leq - \ln(l) )$ into $\nu$ is completely independent of $\tau$. In particular, \eqref{eq:LMExpl} yields that $\tau$ is the left-monotone embedding of $\mu$ into $\nu$ if and only if $(B_0,B_{\tau})$ is the left-curtain coupling of $\mu$ and $\nu$ (cf.\ \cite{BeJu16}). 
	
	The shadow $\shadow{\nu}{\eta}$ of a measure $\eta$ in the probability measure $\nu$, is  the most concentrated (in the sense of $\leqc$) submeasure of $\nu$ which can be reached by an embedding of $\eta$ into $\nu$ via a (randomized) $\mathcal{F}$-stopping time (cf.\ Lemma \ref{lemma:ConvOrder}). Hence, Theorem \ref{thm:intro} characterizes in general barrier solutions as those solutions $\tau$, for which there exists a random time given by $(T_l)_{l \geq 0}$ such that for all $l \geq 0$ the mass which is not stopped before $T_l$ under $\tau$, is allocated by $\tau$ as concentrated as possible in the target distribution $\nu$ without interference with the mass that is stopped before $T_l$.
	
	\subsection*{Interpolation}
	
	If the time-change  $(T_l)_{l \geq 0}$ is measurable w.r.t.\ the completion of the natural filtration $\mathcal{F}^B$ generated by the Brownian motion (as it is the case for the Root-embedding and the left-monotone embedding), we can assume that the Brownian motion $B$ is defined on the canonical path space $\Omega = C([0, \infty))$ and we can consider the natural shift operator $\theta$. In this case, for all $\lambda \in (0,\infty)$ we obtain an interpolation $(R^\lambda_l)_{l \geq 0}$ between two $\mathcal{F}$-time-changes $(T_l ^1)_{l \geq 0}$ and $(T_l^2)_{l \geq 0}$ by
	\begin{equation*}
	R^\lambda_l := T^1 _{l \land \lambda} + (T^2_{l-\lambda} \circ \theta_{T^1_\lambda}) \1_{\{l \geq \lambda\}} = \begin{cases}
	T_l^1 & l \leq \lambda \\
	T_\lambda^1 + T^2 _l \circ  \theta_{T^1_\lambda} & l > \lambda
	\end{cases}.
	\end{equation*}
	
	For the Root time-change $(T_l ^{r})_{l \geq 0}$ and the left-monotone time-change $(T_l^{lm})_{l \geq 0}$  the interpolation becomes
	\begin{equation*}
	R^{\lambda} _l :=
	T^{r}_{l \land \lambda} + (T^{lm} _{l - \lambda} \circ \theta _{T^{r}_{\lambda}}) \1_{\{l \geq \lambda\}} = \begin{cases}
	l &  l \leq \lambda \\
	l & \exp(-B_\lambda) + \lambda \geq l > \lambda \\
	+ \infty  & \exp(-B_\lambda) + \lambda < l, l > \lambda
	\end{cases}.
	\end{equation*}
	A solution $\tau ^{\lambda}$ to $\SEP(\mu,\nu)$ that satisfies property \eqref{eq:ShadowResid} w.r.t.\ $(R^\lambda _l)_{l \geq 0}$, is by Theorem \ref{thm:intro} a barrier solution w.r.t.\ the level-process
	\begin{equation} \label{eq:lvlPrc}
	X_t ^\lambda := \sup \{l \geq 0 : R_l^\lambda \leq t\} = \begin{cases}	
	t & t < \lambda \\
	\lambda + \exp(-B_0)   & t \geq \lambda	\end{cases}.
	\end{equation} 
	A natural guess is that $\lambda \mapsto \tau ^{\lambda}$ is a reasonable interpolation between the left-monotone embedding ($\lambda \uparrow + \infty$) and the Root embedding $(\lambda \downarrow 0)$. This is indeed the case:

	\begin{proposition} \label{prop:Interpolation}
		Let $\lambda \in (0, \infty)$. We define the stochastic process $(X^{\lambda}_t)_{t \geq 0}$ as in \eqref{eq:lvlPrc}.
		There exists a barrier $\mathcal{R}^{\lambda} \subset [0, \infty) \times \mathbb{R}$ such that the first hitting time 
		\begin{equation*}
		\tau ^{\lambda} := \inf \{t \geq 0: (X_t ^{\lambda}, B_t) \in \mathcal{R}^{\lambda} \}
		\end{equation*}
		is a solution to $\SEP(\mu,\nu)$. Moreover, $\Law(B,\tau ^\lambda)$  (as a measure on $\Omega \times [0, \infty))$, converges weakly to $\Law(B,\tau ^r)$ as $\lambda \rightarrow \infty$ and, if $\mu$ is atomless, converges weakly to $\Law(B,\tau ^{lm})$ as $\lambda \rightarrow 0$.
	\end{proposition}
	
	\begin{figure} 
		\begin{center}
			\begin{tikzpicture}
			\draw[->] (0,-2.1) -- (0,2.1);
			\node[above] at (0,2) {$B_t$};
			\draw[->] (0,0) -- (5,0);
			\node[right] at (5,0) {$X ^{\lambda_1}_t$};
			\draw[dashed,thin] (1.5,-2) -- (1.5,2);
			\draw[dotted, domain = 1.75:5,smooth,variable=\t] 
			plot({\t},{-log2(\t-1.5)});
			\fill[fill = cyan, opacity = 0.2] (0,2) -- (5,2) -- (5,0) -- (4,0) -- (2,1) -- (0,1.8);
			\fill[fill = cyan, opacity = 0.2] (0,-2) -- (5,-2) -- (5,0) -- (4,0) -- (4,-0.5) -- (3,-1.25) -- (0,-1.8);
			\node[below] at (1.5,-2) {$\lambda_1$};
			\node[below] at (4,2) {$R^{\lambda_1}$};

			\draw[red] (0,0) -- ++ (0.1,0.2) -- ++ (0.1,-0.1) -- ++ (0.1,0.1) -- ++ (0.1,-0.3) -- ++ (0.1,0.4) -- ++ (0.1,-0.1) -- ++ (0.1,0.4) -- ++ (0.1,-0.1) -- ++ (0.1,-0.1) -- ++ (0.1,0.3) -- ++ (0.1,-0.5)-- ++ (0.1,0.2)-- ++ (0.1,-0.3)-- ++ (0.1,-0.1) -- ++ (0.1,0.2);
			
			\draw[green] (0,-0.5) -- ++ (0.1,-0.2) -- ++ (0.1,0.1) -- ++ (0.1,-0.1) -- ++ (0.1,0.2) -- ++ (0.1,0.1) -- ++ (0.1,-0.4) -- ++ (0.1,0.1) -- ++ (0.1,-0.5) -- ++ (0.1,-0.1) -- ++ (0.1,0.3) -- ++ (0.1,-0.5)-- ++ (0.1,0.2)-- ++ (0.1,-0.3);

			\draw[red] (2.305,-1) -- (2.305, 0.85);
			\end{tikzpicture}
			\begin{tikzpicture}
			\draw[->] (0,-2.1) -- (0,2.1);
			\node[above] at (0,2) {$B_t$};
			\draw[->] (0,0) -- (5,0);
			\node[right] at (5,0) {$X ^\lambda_t$};
			\draw[dashed,thin] (2.5,-2) -- (2.5,2);
			\draw[dotted, domain = 2.75:5,smooth,variable=\t] 
			plot({\t},{-log2(\t-2.5)});
			\fill[fill = cyan, opacity = 0.2] (0,2) -- (5,2) -- (5,0) -- (3.5,0) -- (3,0.3) -- (2,0.9)  -- (1.5,1.2)-- (0,1.8);
			\fill[fill = cyan, opacity = 0.2] (0,-2) -- (5,-2) -- (5,0) -- (3.5,0) -- (4,-0.25) -- (4.5, -1) --  (3,-1.25) -- (0,-1.8);
			\node[below] at (2.5,-2) {$\lambda_2$};
			\node[below] at (4,2) {$R^{\lambda_2}$};

			\draw[red] (0,0) -- ++ (0.1,0.2) -- ++ (0.1,-0.1) -- ++ (0.1,0.1) -- ++ (0.1,-0.3) -- ++ (0.1,0.4) -- ++ (0.1,-0.1) -- ++ (0.1,0.4) -- ++ (0.1,-0.1) -- ++ (0.1,-0.1) -- ++ (0.1,0.3) -- ++ (0.1,-0.5)-- ++ (0.1,0.2)-- ++ (0.1,-0.3)-- ++ (0.1,-0.1) -- ++ (0.1,0.2) -- ++ (0.1,-0.05) -- ++ (0.1,0.3) -- ++ (0.1,0.1) -- ++ (0.1,-0.05) -- ++ (0.1,0.3) -- ++ (0.05,0.1) ;
			
			\draw[green] (0,-0.5) -- ++ (0.1,-0.2) -- ++ (0.1,0.1) -- ++ (0.1,-0.1) -- ++ (0.1,0.2) -- ++ (0.1,0.1) -- ++ (0.1,-0.4) -- ++ (0.1,0.1) -- ++ (0.1,-0.5) -- ++ (0.1,-0.1) -- ++ (0.1,0.3) -- ++ (0.1,-0.5)-- ++ (0.1,0.2)-- ++ (0.1,-0.3);

			\
			\end{tikzpicture}
			\caption{The sketch of two sample paths of $(X^\lambda_t,B_t)_{t \in [0,\tau ^{\lambda}]}$ in the context of Proposition \ref{prop:Interpolation} for two different $\lambda_1 < \lambda _2$ in $(0, \infty)$.}
		\end{center}
	\end{figure}
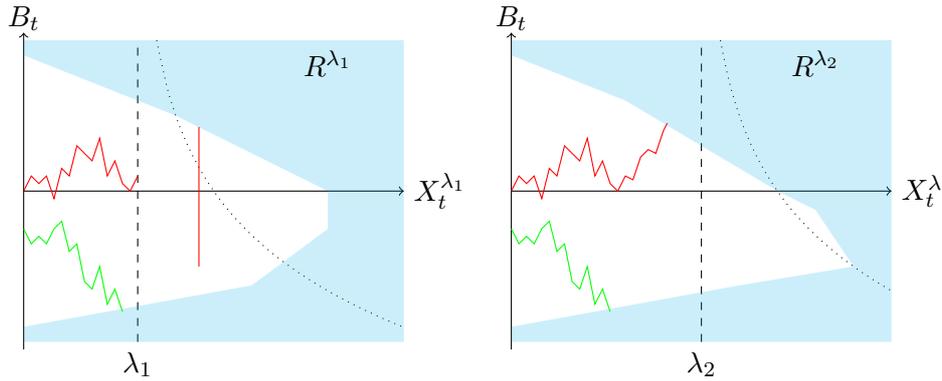
	
	\begin{remark}
		The choice of the Root embedding and the left-monotone embedding as the endpoints of the interpolation is partially arbitrary. As long as both time-changes are $\mathcal{F}^B$-measurable, this procedure can be applied to any two barrier solutions to obtain a new mixed barrier solution (see Lemma \ref{lemma:Nesting}). The continuity and convergence is then a question of the stability properties of the corresponding embeddings.
		
		Other approaches to interpolate (in some sense) between two different barrier solutions can be found in \cite{CoHo07} and \cite{GaObZo19}.
	\end{remark}

	\subsection*{Multi-Marginal Embeddings}
	
	Theorem \ref{thm:intro} can be extended to the case that the barrier solution is ``delayed'', in the sense that the solution can be written as the first hitting time of a barrier after it surpassed a fixed stopping time $\sigma$.
	
	\begin{proposition} \label{prop:ShiftedThm}
		Let $\tau$ be a $\mathcal{F}$-stopping-time that solves $\SEP(\mu,\nu)$. Let $\sigma \leq \tau$ be another $\mathcal{F}$-stopping time. The following are equivalent:
		\begin{itemize}
			\item [(i)] There exists a right-continuous $\mathcal{F}$-adapted stochastic process $(X_t)_{t \geq 0}$ which is non-negative, monotonously increasing and satisfies $\mathbb{P}[\exists s < t : X_s = X_{t} = l] = 0$ for all $l \geq 0$, and a closed barrier $\mathcal{R} \subset [0, \infty) \times \mathbb{R}$ such that 
			\begin{equation*}
			\tau = \inf \{t \geq \sigma : (X_t,B_t) \in \mathcal{R} \} \quad a.s.
			\end{equation*}
			\item [(ii)] There exists a left-continuous $\mathcal{F}$-time-change $(T_l)_{l \geq 0}$ with $T_0 = 0$, $T_\infty = \infty$ and $\mathbb{P}[\lim _{k \downarrow l} T_k = T_l] = 1$ for all $l \geq 0$ such that for all $l \geq 0$ we have
			\begin{equation*}
			\Law(B_\tau; \tau \geq \sigma \lor T_l) = \shadow{\nu}{\Law(B_{\sigma \lor T_l}; \tau \geq \sigma \lor T_l)}.
			\end{equation*}
		\end{itemize}
	\end{proposition}
	
	Motivated by financial applications, there has been an increased interest in the multi-marginal Skorokhod Embedding Problem  and in particular in multi-marginal barrier solutions (cf.\ \cite{BeCoHu17b, NuStTa17}). Since this is essentially a sequence of delayed barrier solutions, we can extend Theorem \ref{thm:intro}  to this case by an inductive application of Proposition \ref{prop:ShiftedThm}.
	
	\begin{corollary} \label{cor:MultiMarginal}
		Let $\mu \leqc \nu _1 \leqc ... \leqc \nu_n$ be greater than $\mu$ in convex order and $\tau _1 \leq ... \leq \tau _n$ an increasing sequence of uniformly integrable $\mathcal{F}$-stopping times such that $\tau _i$ is a solution to $\SEP(\mu, \nu_i)$ for all $1 \leq i \leq n$. The following are equivalent:
		\begin{itemize}
			\item [(i)] There exists a suitable process $(X_t)_{t \geq 0}$, and closed barriers $\mathcal{R}^1,...,\mathcal{R}^n \subset [0, \infty) \times \mathbb{R}$ such that
			\begin{align*}
			\tau ^1 &= \inf\{ t \geq 0: (X_t,B_t) \in \mathcal{R}^1\} \quad\text{ and} \\
			\tau ^i &= \inf\{ t \geq \tau ^{i-1}: (X_t,B_t) \in \mathcal{R}^i\} \quad \text{for all } 1 \leq i \leq n.
			\end{align*}
			\item [(ii)] There exists a suitable time-change $(T_l)_{l \geq 0}$, such that for all $l\geq 0$ we have
			\begin{align*}
			\Law(B_{\tau ^1}; \tau ^1 \geq T_l) &= \shadow{\nu _1}{\Law(B_{T_l}; \tau^1 \geq T_l)} \quad\text{ and} \\
			\Law(B_{\tau ^i}; \tau ^i \geq \tau^{i-1} \lor T_l) &= \shadow{\nu _i}{\Law(B_{\tau ^{i-1} \lor T_l}; \tau ^i \geq \tau^{i-1} \lor T_l)} \quad \text{for all } 1 \leq i \leq n.
			\end{align*}
		\end{itemize}
	\end{corollary}
	
\subsection*{Another Perspective on Theorem \ref{thm:intro}} 

We will prove Theorem \ref{thm:intro} in Section \ref{sec:MainResult} using potential theory. However, there is an alternative point of view on this theorem using Choquet-type representations of the barrier stopping time $\tau$ and the terminal law $\mathsf{Law}(B_\tau)$ of the stopped process. 

The most primitive version of a barrier embedding is a first hitting time of the form 
\begin{equation*}
\tau^F := \inf \{t \geq 0 : B_t \in F\} = \inf \{ t \geq 0 : (t,B_t) \in [0, \infty) \times F\}
\end{equation*}
where $F \subset \mathbb{R}$ is a closed set. The terminal distribution $\Law(B_{\tau ^F})$ w.r.t.\ this stopping time can be characterized using the notion of Kellerer dilations. Given a closed set $F \subset \mathbb{R}$ the Keller dilation is defined as the probability kernel 
\begin{equation} \label{eq:KellererDilation}
K^F(x,dy) = \begin{cases}
\frac{x^+ - x}{x^+ - x^-} \delta x_- + \frac{x - x^-}{x^+ - x^-} \delta_{x^-} \quad & x \not \in F \\
\delta_x & x \in F
\end{cases} 
\end{equation}
where $x^+ = \inf (F \cap [x, \infty])$ and $x^- = \sup (F \cap (-\infty,x])$. As a direct consequence of \cite[Satz 25]{Ke73}, for every closed set $F \subset \mathbb{R}$ a stopping time $\tau$ satisfies $\tau = \tau ^F$ a.e.\ if and only if $\Law(B_\tau) = \Law(B_0)K^F$.  

The main idea behind Theorem \ref{thm:intro}  is now the following: In the same way that a barrier solution $\tau$ can be represented as a composition of first hitting times $(\tau ^{F_t})_{t \geq 0}$ for an increasing family of closed sets $(F_t)_{t \geq 0}$, the terminal  law $\mathsf{Law}(B_\tau)$ w.r.t.\ a  stopping time $\tau$ satisfying the shadow relation \eqref{eq:ShadowResid} can be represented using Kellerer dilations $(K^{F_a})_{a \in [0,1]}$ for an increasing family of closed sets $(F_a)_{a \in [0,1]}$. Since for fixed $F$, $\tau ^F$ and $K^F$ are in a one-to-one correspondence, these two representation -one on the level of stopping times and one on the level of target distributions-  are two sides of the same coin. In fact, up to reparametrization of the index set, these two families can be chosen identical. Let us explain these two representations in more detail.

To keep the notation simple, we will only consider the case of the Root-embedding ($X_t^r = t, T_l ^r = l$). For all $\mathcal{F}^B$-stopping times $\tau_1, \tau _2$ and $s \geq 0$, we define the composition
\begin{equation*}
C_{s}(\tau _1, \tau _2) := \tau _1 \land s + \tau _2 \circ \theta_{\tau _1 \land s}
\end{equation*}
where $\theta$ denotes the shift operator on the path space.
The composition $C_{s}(\tau _1, \tau _2)$ is again a stopping time. We also inductively define the stopping times
\begin{equation*}
C_{s_1, ... , s_n}(\tau_1, ... , \tau _n) := C_{s_n}(C_{s_1, ... , s_{n-1}}(\tau _1, ... , \tau _{n-1}), \tau _n).
\end{equation*}
for $0 \leq s_1 \leq ... \leq s_n$. For all $s > 0$ and for all closed sets $F$, we have $C_s(\tau^F,\tau^F) = \tau^F$. Conversely, if there exists $s \geq 0$ and stopping times $\tau_1, \tau_2$ s.t.\  $\tau^F = C_s(\tau_1,\tau_2)$ for a closed set $F \subset \mathbb{R}$, then $\tau_1 \land s = \tau ^F \land s$ and $\tau _2 \circ \theta_{\tau _1 \land s} = \tau^F \circ \theta_{\tau _1 \land s}$. Therefore, stopping times of the form $\tau^F$ are ``extremal'' or ``atomic'' w.r.t.\ the composition operation $C$.

\begin{lemma}
	Let $\tau$ be a stopping time. The following are equivalent:
	\begin{itemize}
		\item[(i)] There exists a right-barrier $\mathcal{R} \subset [0,\infty) \times \mathbb{R}$ s.t.\ $\tau = \inf \{t \geq 0 : (t,B_t) \in \mathcal{R}\}$. 
		\item[(ii)] There exists an increasing family of closed sets $(F_t)_{t \geq 0}$ such that
		\begin{equation*}
		\tau = \lim _{n \rightarrow \infty} C_{2^{-n}, ... , n} (\tau ^{F_{2^{-n}}}, ... , \tau ^{F_n}).
		\end{equation*}
	\end{itemize}
	In this case a possible right-barrier is given by $\mathcal{R} := \overline{\bigcup _{t \geq 0} [t, \infty) \times F_t}$.
\end{lemma}

The  proof of this equivalence is straightforward using the continuity of Brownian motion. We omit the details.

On the level of measures, we obtain a similar representation of the shadow.
For two probability measures $\zeta_1, \zeta _2$ and all $\alpha \in [0,1]$ the convex combination $(1- \alpha) \zeta_1 + \alpha \zeta_2$ is again a probability measure. By a result of Kellerer \cite[Theorem 1]{Ke73}, for every probability measure $\eta$ the extremal elements of the convex set $\{ \zeta : \eta \leqc \zeta\}$ are given by $\left\{ \eta K^F \, : \, F \subset \mathbb{R} \text{ closed}\right\}$.

\begin{lemma}
	Let $\tau$ be a stopping time and set $l_b := \sup \{l \geq 0: \mathbb{P}[\tau \geq l] \geq b\}$ for $b \in [0,1]$. The following are equivalent:
	\begin{itemize}
		\item [(i)] For all $l \geq 0$ we have $\Law(B_\tau; \tau \geq l) = \shadow{\nu}{\Law(B_l; \tau \geq l)}$.
		
		\item [(ii)] There exists an increasing family of closed sets $(F_a)_{a \in [0,1]}$  such that
		\begin{equation*}
		\Law(B_\tau) = \int _0 ^1 \eta_{1-a} K^{F_a} \de a
		\end{equation*}
		where the probability measures $\eta_a$ are defined by $\eta _a := \lim_{\varepsilon \rightarrow 0} \varepsilon ^{-1} \left( \overline{\eta}^{a + \varepsilon} - \overline{\eta}^a \right)$, $a \in [0,1]$, and $\overline{\eta}^\alpha := \Law(B_\tau; \tau \geq l_\alpha) - \frac{\mathbb{P}[\tau \geq l_\alpha] - \alpha}{\mathbb{P}[\tau = l_a]}\Law(B_\tau; \tau = l_\alpha)$ for $\alpha \in [0,1]$.
	\end{itemize}
	In this case we have $\shadow{\nu}{\Law(B_l; \tau \geq l_b)} = \int_0 ^b \eta_{1-a} K^{F_a} \de a$ for all $b \in [0,1]$.
\end{lemma}

Similar to \cite[Proposition 2.7]{BeJu16b} one can show that (i) implies (ii). The reversed implication is an application of Lemma \ref{lemma:ShadowDecomp}. We leave the details to the reader.

	\section{Related Literature}
	
	The Skorokhod Embedding Problem goes back to Skorokhod's work \cite{Sk65} in 1965. After his own solution to the embedding problem, this problem gained considerable attention in the literature and a wide range of different embeddings exploiting different mathematical tools were found. The survey \cite{Ob04} alone covers more than 20 different solutions. Moreover, several interesting variants of the Skorokhod Embedding are considered. Recently, there is an increased interest in a variant of Skorokhod Embedding Problem, which asks embeddings to minimize or maximize a predetermined cost function of space and time. This variant of the Skorokhod Embedding Problem has a direct connection to robust mathematical finance which was first noticed by Hobson \cite{Ho98}. For further background we refer to \cite{Ho03}. A novel mathematical exploration of properties of the optimal Skorokhod Embedding Problem in combination with optimal transport can be found in \cite{BeCoHu17}. Further variants are for instance the extensions to the embedding of multiple distributions (cf.\ \cite{BeCoHu17b}) and to higher dimensions (cf.\ \cite{GhKiPa19}).
	
	Among the first solutions to the Skorokhod Embedding Problem was Root's construction \cite{Ro69} of a barrier solution in the time-space phase space in 1969. Shortly after, Rost \cite{Ro76} proved that the Root-embedding is the  unique embedding which has minimal variance among all other embeddings and provided an alternative construction of this embedding based on the potential theory for Markov processes. The Root-embedding and properties of the corresponding barrier are still subject of current research \cite{CoWa13,GaObZo19}. Moreover, the Root-embedding was recently used to construct a counterexample to the Cantelli-conjecture \cite{KlKu15}. The Root-embedding is presumably the most prominent barrier solution to the Skorokhod Embedding Problem. However, there are several other embeddings which can be characterized as first hitting times of barriers in a different phase space \cite{AzYo79, Va83}.

	The shadow for finite measures on the real line was introduced by Beiglböck and Juillet \cite{BeJu16} as the main tool in their construction of the left-curtain coupling. Thereby, they showed important properties as the associativity law and continuity, and coined the name shadow. Nevertheless, the essential concept of the shadow as well as its existence in a very broad framework already appeared in \cite{Ro71}. 
	The shadow is used to study properties of the left-curtain coupling  (cf.\  \cite{BeJu16}, \cite{Ju14},  \cite{HoNo17}, \cite{HoNo21}).
	Furthermore, the shadow can be used to construct and characterize a whole family of martingale couplings on the real line \cite{BeJu16b}, as well as finite-step martingales \cite{NuStTa17} and solutions to the peacock problem \cite{BrJuHu20}. To the best of our knowledge, the only known connection with the Skorokhod Embedding Problem so far is implicitly through the left-monotone embedding because it is uniquely characterized by the property that the induced martingale coupling between the initial and the terminal marginal distribution is precisely the left-curtain coupling.

	\section{Preliminary Results}
	
	\subsection{Notation} \label{ssec:Notation}

	$\Omega$ is a Polish space equipped with the Borel $\sigma$-algebra, $\mathcal{F}$ is a right-continuous filtration on $\Omega$ and $B$ is a $\mathcal{F}$-Bownian motion on the complete filtered probability space $(\Omega, \mathcal{B}(\Omega), \mathbb{P}, \mathcal{F})$. 	We use the notation $\Law(X;A)$ for the (sub-)probability measure which is given by the push-forward of the random variable $X$ under the restriction of $\mathbb{P}$ to the Borel set $A$. Alternatively, we sometimes use the notation $X_{\#}(\mathbb{P}_{|A})$ for this object.
	
	Further, we denote the set of finite (resp.\ probability) measures on a measurable space $\mathsf{X}$ by $\mathcal{M}(\mathsf{X})$ (resp.\ $\mathcal{P}(\mathsf{X})$). In the case $\mathsf{X} = \mathbb{R}$, we denote by $\MO(\mathbb{R})$ (resp.\ $\PO(\mathbb{R})$) the subset of finite (resp.\ probability) measures with finite first moment. 
	We equip $\MO(\mathbb{R})$ with the initial topology generated by the functionals $(I_f)_{f \in C_b(\mathbb{R}) \cup \{\vert \cdot \vert\}}$ where
	\begin{equation*}
	I_f : \MO(\mathbb{R}) \ni \pi \mapsto \int _{\mathbb{R}} f \de \pi \in \mathbb{R},
	\end{equation*}
	$C_b(\mathbb{R})$ is the set of continuous and bounded functions, and $\vert \cdot \vert$ denotes the absolute value function.
	We denote this topology on $\MO(\mathbb{R})$ by $\TO$. 
	
	Finally, we define two order relations on $\MO(\mathbb{R})$. We say  that $\mu  \in \MO(\mathbb{R})$ is smaller than or equal to $\mu ' \in \MO(\mathbb{R})$ in convex order, $\mu \leqc \mu'$, if
	\begin{equation} \label{eq:OrderRelation}
	\int _{\mathbb{R}} \varphi \de \mu \leq \int _{\mathbb{R}} \varphi \de \mu '
	\end{equation} 
	holds for all convex $\varphi$ and $\mu$ is smaller than or equal to $\mu'$ in positive order, $\mu \leqp \nu$, if \eqref{eq:OrderRelation} holds for all non-negative $\varphi$.

	\subsection{Randomized Stopping Times}
	
	The product space $\Omega \times [0,\infty)$ equipped with the product topology and Borel $\sigma$-algebra is again a Polish space.
	
	\begin{definition}
		A randomized stopping time (RST) w.r.t.\ $\mathbb{P}$ is a subprobability measure $\xi$ on $\Omega \times [0, \infty)$ such that the projection of $\xi$ onto $\Omega$ is $\mathbb{P}$ and there exists a disintegration $(\xi_\omega)_{\omega \in \Omega}$ of $\xi$ w.r.t.\ $\mathbb{P}$ such that 
		\begin{equation} \label{eq:DecompRST}
		\rho _u : \omega \mapsto \inf\{t \geq 0 : \xi_\omega[0,t] \geq u\}
		\end{equation}
		is an $\mathcal{F}$-stopping time for all $u \in [0,1]$. We call a RST $\xi$ finite, if $\xi$ is a probability measure.
	\end{definition}
	
	We equip the space of RST with the topology of weak convergence of measures on $\Omega \times [0, \infty)$, i.e.\ the continuity of functionals $\xi \mapsto \int \varphi \de \xi$ for all $\varphi \in C_b(\Omega \times [0, \infty))$. The RST-property is closed under this topology (cf.\ \cite[Corollary 3.10]{BeCoHu17}).

	Any $\mathcal{F}$-stopping time $\tau$ naturally induces a RST by $\xi^{\tau} := \Law_{\mathbb{P}}(B,\tau)$. Conversely, we can represent any randomized stopping time as a usual stopping time by enlarging the filtration.
	
	\begin{lemma} [{\cite[Theorem 3.8]{BeCoHu17}}] \label{lemma:ReprRST}
		For every RST $\xi$ there exists an $(\mathcal{B}([0,1]) \times \mathcal{F}_t)_{t \geq 0}$-stopping-time $\overline{\tau} ^\xi$ on the probability space $([0,1] \times \Omega, \mathcal{B}([0,1] \times \Omega), \overline{\mathbb{P}})$ where $\overline{\mathbb{P}}$ is the product of the Lebesque measure and $\mathbb{P}$ such that 
		\begin{equation*}
		\xi = \Law_{\overline{\mathbb{P}}}(\overline{\mathsf{Id}},\overline{\tau}^\xi)
		\end{equation*}
		where $\overline{\mathsf{Id}} : (u,\omega) \mapsto \omega$. Moreover, $\overline{B} :  (u, \omega) \mapsto B(\omega)$ is a Brownian motion on $([0,1] \times \Omega, \mathcal{B}([0,1] \times \Omega), \overline{\mathbb{P}})$.
	\end{lemma}
	
	This representation is useful to justify the application of known theorems of stopping times to RST and will be used in the following. For further literature on randomized stopping times we refer to \cite{BeCoHu17} and references therein.

	Provided that $\Law(B_0) = \mu \leqc \nu$, we say that $\xi$ is a solution of $\SEP(\mu,\nu)$ if
	\begin{align*}
	\sup _{s \geq 0} \int _{\Omega \times [0,\infty)} B_{s \land t}  \de \xi(\omega,t) < + \infty \quad \text{and} \quad
	((\omega,t) \mapsto B_t(\omega))_{\#} \xi = \nu.
	\end{align*}
	If $\xi$ is induced by a $\mathcal{F}$-stopping time $\tau$, this definition is consistent with the definiton of $\SEP(\mu,\nu)$ in the introduction. 
	Especially in Section \ref{sec:MainResult} we will use the notational convention that $(\omega,t)$ always refers to an element of $\Omega \times [0, \infty)$. In particular, we will write $\xi[t \geq X]$ instead of $\xi[\{(\omega,t) : t \geq X(\omega)\}]$ where $X$ is a random variable and $\xi$ a RST.

	\subsection{Potential Theory} \label{ssec:PotentialTheory}
	
	Potential Theory is known to be a useful tool when dealing with barrier solutions (cf.\ \cite{Ro76}, \cite{CoWa13}) and the shadow (cf.\ \cite{Ro71}, \cite{BeJu16}). Since it is also a central part of our proof of Theorem \ref{thm:intro}, we recall some results below.
	
	\begin{definition} \label{def:PotentialFunction}
		Let $\eta \in \MO$. The potential function of $\eta$ is defined by
		\begin{equation*}
		U_\eta : \mathbb{R} \rightarrow [0, \infty) \quad U_\eta (x) := \int _{\mathbb{R}} |y - x| \de \eta (y).
		\end{equation*} 
	\end{definition}
	
	Since elements of $\MO$ have finite first moments, the potential function is always well-defined.

	\begin{lemma}[{cf.\ \cite[Proposition 4.2]{BeJu16}, \cite[p.\ 335]{Ob04}}] \label{lemma:ConvOrder}
		Let  $\mu, \nu \in \mathcal{P}_1(\mathbb{R})$. The following are equivalent:
		\begin{itemize}
			\item [(i)] $\mu \leqc \nu$
			\item [(ii)] $U_\mu \leq U_\nu$
			\item [(iii)] There exists a solution to $\SEP(\mu,\nu)$.
		\end{itemize}
	\end{lemma}
	
	The equivalence between (i) and (ii) is not restricted to probability measures. Since both the convex order and the order of the potetntial functions are invariant w.r.t.\ scaling with positive factors, for all $\eta, \zeta \in \MO$ with $\eta(\mathbb{R}) = \zeta(\mathbb{R})$ we have $\eta \leqc \zeta$ if and only if $U_{\eta} \leq U_\zeta$ .

	\begin{lemma}[{cf.\ \cite[Proposition 4.1]{BeJu16}}]
		\label{lemma:characPotF}
		Let $m \in [0,\infty)$ and $x^* \in \mathbb{R}$.
		For a function $u:\mathbb{R} \rightarrow \mathbb{R}$ the following statements are equivalent:
		\begin{enumerate}
			\item [(i)] There exists a finite measure $\mu \in \MO$ with mass $\mu(\mathbb{R}) = m$ and barycenter $x^* = \int _{\mathbb{R}} x \de \mu (x)$ such that $U_\mu = u$ .
			\item [(ii)] The function $u$ is non-negative, convex and satisfies
			\begin{equation} \label{eq:characPotF}
			\lim _{x \rightarrow \pm \infty} u(x) - m|x - x^*| = 0. 
			\end{equation}
		\end{enumerate}
		Moreover, for all $\mu, \mu' \in \MO$ we have $\mu = \mu'$ if and only if $U_\mu = U_{\mu'}$.
	\end{lemma}

	\begin{lemma} \label{lemma:PropPotf}
		Let $\eta$ be a positive measure on $\mathbb{R}$. If there exists an $\varepsilon > 0$ such that $U_{\eta}$ is affine on $[x-\varepsilon, x+ \varepsilon]$, $x \not \in \mathrm{supp}(\eta)$.
	\end{lemma}
	
	\begin{proof}
		The claim follows from the observation that the potential function of the measure $\eta$ satisfies $\frac{1}{2} U_\eta '' = \eta$ in a distributional sense (cf.\ \cite[Proposition 2.1]{HiRo12}). 
	\end{proof}
	
	\begin{corollary} \label{cor:EquaPotfToZero}
		Let $\mu \leq \nu$ and $\tau$ be a solution to $\SEP(\mu,\nu)$. We have
		\begin{equation*}
		\mathbb{P}[\tau > 0, U_{\mu}(B_0) = U_{\nu}(B_0)] = 0.
		\end{equation*}
	\end{corollary}
	
	\begin{proof}
		Let $A := \{x \in \mathbb{R} :U_{\mu}(x) = U_{\nu}(x)\}$ and set $\eta := \mathrm{Law}(B_{0}; B_{0} \in A)$. Fubini's Theorem yields 
		\begin{align*}
		0 = \int U_{\nu} - U_{\mu} \de \eta  = \mathbb{E}[U_{\eta}(B_\tau) - U_{\eta}(B_0)] =  \mathbb{E}\left[\left(U_{\eta}(B_\tau) - U_{\eta}(B_0)\right) \1_{\{\tau > 0\}}\right].
		\end{align*}	
		Since $U_\eta$ is a convex function and $(B_{t \land \tau})_{t \geq 0}$ is a uniformly integrable martingale, the (conditional) Jensen inequality yields that 
		$U_\eta$ is $\mathbb{P}$-a.s.\ affine at $B_0$ on the set $\tau > 0$.
		Hence, by Lemma \ref{lemma:PropPotf} the claim follows.
	\end{proof}

	\begin{lemma} \label{lemma:T1Conv}
		Let $(\mu_n)_{n \in \mathbb{N}}$ be a sequence in $\MO(\mathbb{R})$. The following are equivalent:
		\begin{itemize}
			\item [(i)] The sequence $(\mu_n)_{n \in \mathbb{N}}$ is weakly convergent and there exists a finite measure $\eta \in \MO(\mathbb{R})$ such that 
			\begin{equation*}
			\int _{\mathbb{R}} \varphi \de \mu_n \leq \int _{\mathbb{R}} \varphi \de \eta
			\end{equation*}
			for all non-negative convex $\varphi$.
			
			\item [(ii)] The sequence $(\mu_n)_{n \in \mathbb{N}}$ is convergent under $\TO$.
			
			\item [(iii)] The sequence of potential functions is pointwise convergent and the limit is the potential function of a finte measure.
		\end{itemize}
	\end{lemma}
	
	\begin{proof}
		For the equiavlence of (ii) and (iii) and the implication (i)$\Rightarrow$(ii)  we refer to \cite[Lemma 3.6]{BrJuHu20} and \cite[Lemma 3.3]{BrJuHu20}. It remains to show that (ii) implies (i). Since $\TO$ is by definition stonger than the weak topology, $(\mu_n)_{n \in \mathbb{N}}$ is weakly convergent. Moreover, by \cite[Proposition 7.1.5]{AmGiSa08} the convergence in $\TO$ implies that
		\begin{equation*}
		\limsup _{K \rightarrow \infty} \sup _{n \in \mathbb{N}} \int _{\mathbb{R}} |x| \1_{\{\vert x \vert \geq K\}} \de \mu_n(x) = 0.
		\end{equation*}
		Hence, there exists a sequence $(K_m)_{m \in \mathbb{N}}$ with $K_{m+1} \geq K_m \geq 1$  such that $$\sup _{n \in \mathbb{N}} \int _{\mathbb{R}} |x| \1_{\{\vert x \vert \geq K_m\}} \de \mu_n(x) \leq 2^{-m}$$ for all $m \in \mathbb{N}$. The measure
		\begin{equation*}
		\eta := \sum _{m = 1} ^{\infty} \sup _{n \in \mathbb{N}} \mu_n \left( [-K_m,-K_{m-1}] \cup [K_{m-1},K_m] \right) \left( \delta _{-K_m} + \delta_{K_m} \right)
		\end{equation*}
		is an element of $\mathcal{M}_1(\mathbb{R})$ which satisfies the desired properties.
	\end{proof}

	\subsection{Shadows}
	
	Recall the definition of the shadow in Definition \ref{def:Intro}. As direct consequences of this definition we obtain that 
	\begin{equation*}
	\eta \leqp \nu \Rightarrow \shadow{\nu}{\eta} = \eta  \quad \text{and} \quad \eta \leqc \eta' \Rightarrow \shadow{\nu}{\eta} \leqc \shadow{\nu}{\eta'}. 
	\end{equation*}
	In the following we collect further properties of the shadow.

	\begin{lemma}[{\cite[Theorem 4.8]{BeJu16}}] \label{lemma:ShadowAssz}
		Let  $\eta := \eta _1 + \eta _2 \leqc \nu$, the shadow of $\eta_2$ in $\nu - \shadow{\nu}{\eta_1}$ exists and we have
		\begin{equation*}
		\shadow{\nu}{\eta} = \shadow{\nu}{\eta _1} + \shadow{\nu - \shadow{\nu}{\eta _1}}{\eta _2}.
		\end{equation*}
	\end{lemma}
	
	The statement in Lemma \ref{lemma:ShadowAssz} is the ``associativity law'' for shadows already mentioned in the introduction.
	
	\begin{corollary} \label{lemma:CharShad}
		Let $\mu \leqc \nu$ be probability measures and $A \subset \mathbb{R}$ a Borel set such that $\mu(A) > 0$. 
		If a solution $\tau$ of $\SEP(\mu,\nu)$ satisfies 
		\begin{equation*}
		\forall \tau' \text{ solution of } \SEP(\mu,\nu) \, : \, \Law(B_{\tau}; B_0 \in A) \leqc \Law(B_{\tau'}; B_0 \in A) , 
		\end{equation*}
		we have $\Law(B_\tau; B_0 \in A) = \shadow{\nu}{\mu _{|A}}$.
	\end{corollary}
	
	\begin{proof}
		If $\alpha := \mu(A) = 1$, there is nothing to show because $\shadow{\nu}{\mu_{|A}} = \nu = \Law(B_\tau; B_0 \in A)$. 
		Assume $\alpha < 1$.
		Since $\tau$ is a solution to $\SEP(\mu,\nu)$, we have
		\begin{equation*}
		\mu_{|A} = \Law(B_0; B_0 \in A) \leqc \Law(B_\tau; B_0 \in A) \leqp \nu
		\end{equation*}
		and hence we obtain $\shadow{\nu}{\mu_{|A}} \leqc \Law(B_\tau; B_0 \in A)$. It remains to show that also the reversed relation holds.
		By definition of the shadow, we have $\mu_{|A} \leqc \shadow{\nu}{\mu_{|A}}$ and Lemma \ref{lemma:ConvOrder} yields that there exists a solution $\tau^A$ to $\SEP(\alpha^{-1}\mu_{|A}, \alpha^{-1}\shadow{\nu}{\mu_{|A}})$. By Lemma \ref{lemma:ShadowAssz} it is
		\begin{equation*}
		\mu_{|A^c} \leq_c \nu - \shadow{\nu}{\mu_{|A}}
		\end{equation*}
		and again Lemma \ref{lemma:ConvOrder} yields the existence of a solution $\tau ^{A^c}$ to $\SEP((1-\alpha)^{-1}\mu_{|A^c}, (1-\alpha)^{-1}(\nu - \shadow{\nu}{\mu_{|A}}))$. Since $\{B_0 \in A\} \in \mathcal{F}_0$, 
		\begin{equation*}
		\tau' := \tau ^A \1_{\{B_0 \in A\}} + \tau ^{A^c} \1_{\{B_0 \not \in A\}}
		\end{equation*}
		is a solution to $\SEP(\mu,\nu)$ and thus
		\begin{equation*}
		\Law(B_{\tau}; B_0 \in A) \leqc \Law(B_{\tau'}; B_0 \in A) = \alpha \Law(B_{\tau^A}) = \shadow{\nu}{\mu_{|A}}. \qedhere
		\end{equation*}
	\end{proof}
	
	\begin{corollary} \label{cor:ShadowOnEqualPart}
		Let $\mu \leqc \nu$ and $\tau$ be a solution to $\SEP(\mu,\nu)$. Let $A \in \mathcal{F}_0$ such that $U_{\mu}(B_0) = U_{\nu}(B_0)$ on $A$. Then 
		\begin{equation*}
		\shadow{\nu}{\Law(B_0; A^c)} = \Law(B_\tau; A^c).
		\end{equation*}
	\end{corollary}
	
	\begin{proof}
		Set $I := \{ x \in \mathbb{R} : U_{\mu}(x) < U_{\nu}(x)\}$. Since $I$ is the collection of irreducible components of $(\mu,\nu)$ (cf.\ \cite[Section A.1]{BeJu16} ),  for any solution $\tau'$ of $\mathrm{SEP}(\mu,\nu)$, the stopped process $(B_{\tau' \land s})_{s \geq 0}$ stays in the irreducible component that it started in. Hence, the measure $\Law(B_{\tau'};B_0 \in I)$ is independent of the specific solution $\tau'$.
		By Lemma \ref{lemma:CharShad}, we obtain
		\begin{equation} \label{eq:IrredShadow}
		\mathrm{Law}(B_{\tau'}; B_{0} \in I) = \shadow{\nu}{\mathrm{Law}(B_{0}; B_{0} \in I)}
		\end{equation}
		for any solution $\tau'$ of $\SEP(\mu,\nu)$.
		
		Since  $\{B_0 \in I\} \subset A^c$ and $\tau = 0$ on $\{B_0 \not \in I\}$ by Corollary \ref{cor:EquaPotfToZero}, we obtain 
		\begin{align*}
		\Law(B_0; B_0 \not \in I, A^c) = \Law(B_\tau; B_0 \not \in I, A^c) \leq_+ \Law(B_\tau; B_0 \not \in I). 
		\end{align*}
		Thus, with Lemma \ref{lemma:ShadowAssz} and \eqref{eq:IrredShadow} we obtain
		\begin{align*}
		\shadow{\nu}{\Law(B_0; A^c)} &= \shadow{\nu}{\Law(B_0; B_0 \in I)} + \shadow{\nu - \shadow{\nu}{\Law(B_0; B_0 \in I)}}{\Law(B_0; B_0 \not \in I, A^c)} \\
		&= \Law(B_\tau ; B_0 \in I) + \shadow{\Law(B_\tau; B_0 \not \in I)}{\Law(B_0; B_0 \not \in I, A^c)} \\
		&= \Law(B_\tau ; B_0 \in I) + \Law(B_\tau; B_0 \not \in I, A^c) = \Law(B_\tau ; A^c). \qedhere
		\end{align*}
	\end{proof}
	
	The connection of shadows to potential theory is through the following characterization of the potential functions of the shadow.
	
	\begin{lemma}[{\cite[Theorem 2]{BeHoNo20}}] \label{lemma:PotfShad}
		Let $\hat{\mu} \leq \mu \leqc \nu$. The potential function of the shadow $\shadow{\nu}{\hat{\mu}}$ is given by
		\begin{equation*}
		U_{\shadow{\nu}{\hat{\mu}}} = U_{\nu} - \mathrm{conv} \left( U_{\nu} - U_{\hat{\mu}} \right)
		\end{equation*}
		where $\mathrm{conv}(f)$ denotes the convex hull of a function $f$, i.e. the largest convex function that is pointwise smaller than $f$.
	\end{lemma}

	\begin{lemma} [{\cite[Lemma 1]{BeHoNo20}}] \label{lemma:PropConv}
		Let $f$ be a continuous function bounded by an affine function from below. If $x \in \mathbb{R}$ satisfies $(\mathrm{conv}(f))(x) < f(x)$, there exists an $\varepsilon > 0$ such that $\mathrm{conv}(f)$ is affine on $[x - \varepsilon, x + \varepsilon]$.
	\end{lemma}

	\begin{lemma} \label{lemma:ShadowDecomp} 
		Let $(\mu_a)_{a \in [0,1]}$ be a family of probability measures, $(F_a)_{a \in [0,1]}$ a decreasing sequence of closed subsets of $\mathbb{R}$ and set $\nu = \int _0 ^1 \mu_a K^{F_a} \leqp \nu$. For all $b \in [0,1]$ we have
		\begin{equation*}
		\mathcal{S}^{\nu}\left(\int _0 ^b \mu_a \de a\right) = \int _0 ^b \mu_a K^{F_a} \de a.
		\end{equation*}
	\end{lemma}
	
	\begin{proof}
		Let $\eta, \zeta \in \MO(\mathbb{R})$ and $F \subset \mathbb{R}$ a closed set with $\mathrm{supp}(\zeta) \subset F$. Since we have 
		\begin{equation*}
		\eta \leqc \eta K^F \leqp \eta K^F + \zeta,
		\end{equation*}
		we obtain $\shadow{\eta K^F + \zeta}{\eta} \leqc \eta K^F$. Conversely, we also have 
		\begin{equation*}
		\eta K^F \leqc  \eta K^{\mathrm{supp}(\eta K^F + \zeta)} \leqc \shadow{\eta K^F + \zeta}{\eta}  
		\end{equation*}
		because $\mathrm{supp}(\eta K^F + \zeta) \subset F$ and by definition $\eta K^{\mathrm{supp}(\eta K^F + \zeta)}$ is the smallest measure in convex order which dominates $\eta$ in convex order and is supported on $\mathrm{supp}(\eta K^F + \zeta)$ (cf.\ \eqref{eq:KellererDilation}). Hence, we have $ \shadow{\eta K^F + \zeta}{\eta} = \eta K^F$.
		Furthermore, for all $n \in \mathbb{N}$, $\mu_1, \ldots , \mu_n \in \MO$ and closed sets $F_1, \ldots , F_n \subset \mathbb{R}$ we can apply this equality to get
		\begin{equation*}
		\mu_1 K^{F_1} = \shadow{\mu_1 K ^{F_1} + \ldots + \mu_n K^{F_n}}{\mu _1} 
		\end{equation*}
		and with Lemma \ref{lemma:ShadowAssz} we inductively obtain
		\begin{align*}
		&\shadow{\mu_1 K ^{F_1} + \ldots + \mu_n K^{F_n}}{\mu _1 + \ldots + \mu_{k}} \\
		&= \shadow{\mu_1 K ^{F_1} + \ldots + \mu_n K^{F_n}}{\mu _1 + \ldots + \mu_{k-1}} \\
		& \quad \quad + \shadow{\mu_1 K ^{F_1} + \ldots + \mu_n K^{F_n} - \shadow{\mu_1 K ^{F_1} + \ldots + \mu_n K^{F_n}}{\mu _1 + \ldots + \mu_{k-1}} }{\mu_{k}} \\
		&= \mu _1K^{F_1} + \ldots + \mu_{k-1}K^{F_{k-1}} + \shadow{\mu_{k} K ^{F_{k}} + \ldots + \mu_n K^{F_n}}{\mu_{k}} \\
		&= \mu_1 K ^{F_1} + \ldots + \mu_k K^{F_k}
		\end{align*}
		for all $2 \leq k \leq n$.
		Since the map $(\mu,\nu) \mapsto \shadow{\nu}{\mu}$ is continuous under $\TO$ (cf.\  \cite{Ju14}), the claim follows.
	\end{proof}

	\section{Proof of the Main Result} \label{sec:MainResult}
	
	We split the proof of Theorem \ref{thm:intro} in three parts. In Subsection \ref{ssec:adjoint} we show that the assumptions on the time-change and the level process in Theorem \ref{thm:intro} correspond to each other. In Subsection \ref{ssec:AprioriBound} we construct for every solution of the Skorokhod Embedding Problem an upper bound in the form of a barrier solution and we prove in Subsection \ref{ssec:ActualProof} that this upper bound is attained if and only if the properties of Theorem \ref{thm:intro} are satisfied. 
	
	\subsection{Monotonously Increasing Processes} \label{ssec:adjoint}
	
	\begin{definition}
		Two monotonously increasing and non-negative families of random variables $(X_t)_{t \geq 0}$ and $(T_l)_{l \geq 0}$ are adjoint if $\mathbb{P}[X_t \geq l \Leftrightarrow T_l \leq t] = 1$ for all $l,t \geq 0$.
	\end{definition}
	
	\begin{remark}
		If $(X_t)_{t \geq 0}$ is right-continuous or $(T_l)_{l \geq 0}$ left-continuous and both families are adjoint, we have $\mathbb{P}[\forall l,t \geq 0 : X_t \geq l \Leftrightarrow T_l \leq t] = 1$.
	\end{remark}

	\begin{lemma} \label{lemma:ExAdjoint}
		\begin{itemize}
			\item [(i)] Let $(X_t)_{t \geq 0}$ be a right-continuous $\mathcal{F}$-adapted stochastic process which is non-negative, monotonously increasing and satisfies $\mathbb{P}[\exists s < t : X_s = X_t = l] = 0$ for all $l \geq 0$. Then, the family $(T_l)_{l \geq 0}$ defined by 
			\begin{equation*}
			T_l := \inf \{t \geq 0 : X_t \geq l \}
			\end{equation*}
			is a left-continuous $\mathcal{F}$-time change with $T_0 = 0$, $T_\infty = + \infty$ and $\mathbb{P}[\lim _{k \downarrow l} T_k = T_l] = 1$ for all $l \geq 0$ which is adjoint to $(X_t)_{t \geq 0}$.
			
			\item [(ii)] Let $(T_l)_{l \geq 0}$ be a left-continuous $\mathcal{F}$-time-change  with $T_0 = 0$, $T_\infty = \infty$ and $\mathbb{P}[\lim _{k \downarrow l} T_k = T_l] = 1$ for all $l \geq 0$. Then. the family $(X_t)_{t \geq 0}$ defined by
			\begin{equation*}
			X_t := \sup \{l \geq 0 : T_l \leq t\}
			\end{equation*}
			is a right-continuous $\mathcal{F}$-adapted stochastic process $(X_t)_{t \geq 0}$ which is non-negative, monotonously increasing and satisfies $\mathbb{P}[\exists s < t: X_s = X_{t} = l] = 0$ for all $l \geq 0$, and which is adjoint to $(T_l)_{l \geq 0}$.
		\end{itemize}
	\end{lemma}
	
	\begin{proof}
		Item (i): Let $t,l \geq 0$. If $X_t \geq l$, $T_l \leq t$ directly by definition. Conversely, if $T_l \leq t$, for all $u > t$ we obtain $X_u \geq l$ and thus $X_t = \lim _{u \downarrow t} X_u \geq l$ by right-continuity of $X$. Hence, $(T_l)_{l \geq 0}$ is adjoint to $(X_t)_{t \geq 0}$.
		
		Clearly, $(T_l)_{l \geq 0}$ is monotonously increasing.
		Since $(T_l)_{l \geq 0}$ and $(X_t)_{t \geq 0}$ are adjoint, the symmetric difference $\{T_l \leq t\} \triangle \{X_t \geq l\}$ is a $\mathbb{P}$-null-set and therefore contained in the completed filtration $\mathcal{F}_t$. Thus,  $(T_l)_{l \geq 0}$ is a $\mathcal{F}$-time-change. Since $X_t$ is non-negative and finite, we obtain $T_0 = 0$ and $T_{\infty} = + \infty$. Moreover, $l \mapsto T_l$ is left-continuous by definition. 
		
		Furthermore, we have $\mathbb{P}[ \lim _{k \downarrow l} T_k > T_l ] \leq \mathbb{P}[\exists s < t : X_s = X_{t} = l] = 0$.
		
		Item (ii): Basically the same just in reverse.
	\end{proof}
	
	\textbf{In the following} we fix a $\mathcal{F}$-adapted stochastic process $(X_t)_{t \geq 0}$ and an adjoint $\mathcal{F}$-time-change $(T_l)_{l \geq 0}$ that satisfy the properties listed in Lemma \ref{lemma:ExAdjoint}.

	\subsection{A-priori Bound} \label{ssec:AprioriBound}
	
	Let $B$ be a Brownian motion that starts in $\mu$.
	Fix a randomized stopping time $\xi$ that is a solution to $\SEP(\mu,\nu)$. 
	
	To simplify notation we will use the following notation for measures derived from $\xi$
	\begin{equation} \label{eq:DefRST2}
	\Law(B_{\sigma \land \xi}) := ((\omega,t) \mapsto B_{\sigma(\omega) \land t}(\omega))_{\#} \xi 
	\end{equation}
	where $\sigma$ is an $\mathcal{F}$-stopping time.
	We set $u(l,x) := U_{\Law(B_{T_l \land \xi})}(x)$ and $v(x) := U_{\Law(\nu)}$ for $l \geq 0$ and $x \in \mathbb{R}$.
	In this part we will show that $\xi$ is bounded from above by the stopping time 
	\begin{equation*}
	\hat{\tau} := \inf\{ t \geq 0 : u(X_t,B_t) = v(B_t)\},
	\end{equation*}
	i.e.\ we have $\xi[t \leq \hat{\tau}] = 1$. Since $u$ depends on $\xi$, $\hat{\tau}$ is obviously not a global bound for all solutions to $\mathrm{SEP}(\mu,\nu)$. Nevertheless, Lemma \ref{lemma:uCont} implies that $\hat{\tau}$ is a barrier solution.

	\begin{lemma} \label{lemma:uCont}
		The function $u$ is continuous and monotonously increasing in the first component. Moreover, for all $x \in \mathbb{R}$ we have $v(x) = \lim_{l \rightarrow \infty} u(l,x)$.
	\end{lemma}
	
	\begin{proof}
		For all $x \in \mathbb{R}$ and $l \leq l'$, by Lemma \ref{lemma:ReprRST} we have
		\begin{align*}
		u(l',x) =  \overline{\mathbb{E}} \left[ |\overline{B}_{T_{l'} \land \overline{\tau}^\xi} - x| \right] \geq \overline{\mathbb{E}} \left[ |\overline{B}_{T_{l} \land \overline{\tau}^\xi}  - x| \right]  = u(l,x)
		\end{align*}
		because $\Law_{\overline{\mathbb{P}}}(\overline{B}_{T_{l} \land \overline{\tau}^\xi})_{l \geq 0}$ is increasing in convex order by the optional stopping theorem.
		We chose $(T_l)_{l \geq 0}$ such that, for fixed $l_0 \geq 0$, $l \mapsto T_l$ is $\mathbb{P}$-a.s.\ continuous at $l_0$. Hence, $l \mapsto \Law(B_{T_l \land \xi})$ is weakly continuous and by Lemma \ref{lemma:T1Conv}, $u$ is continuous in the fist component because $\Law(B_{T_l \land \xi}) \leq_c \nu$ for all $l \geq 0$. Furthermore, $u$ is $1$-Lipschitz continuous in the second component because $u(l,\cdot)$ is the potential function of $\Law(B_{T_l \land \xi})$.
	\end{proof}

	\begin{lemma} \label{lemma:EquaToLeq}
		Let $l \geq 0$ and $\sigma$ be a finite $\mathcal{F}$-stopping time. It is
		\begin{equation*}
		\xi \left[  u(l,B_{\sigma}) = v(B_{\sigma}), t > \sigma \geq T_l \right] = 0.
		\end{equation*}
	\end{lemma}
	
	\begin{proof}
		This is a direct consequence of Lemma \ref{lemma:ReprRST} and Corollary \ref{cor:EquaPotfToZero}.
	\end{proof}

	\begin{proposition} \label{prop:EquaToLeq}
		Let $\sigma$ be a finite $\mathcal{F}$ stopping time with $\mathbb{P}[u(X_{\sigma},B_{\sigma}) = v(B_{\sigma})] = 1$. We have $\xi[t \leq \sigma] = 1$.
	\end{proposition}
	
	\begin{proof}
		Let $r(x) := \inf \{l \geq 0 : u(l,x) = v(x)\}$ for all $x \in \mathbb{R}$ and let 
		\begin{equation} \label{eq:DefOfL}
		L := \{ r(x) : x \in \mathbb{R}, \exists \varepsilon > 0 \text{ s.t. } r(x) \leq r(y) \text{ f.a.\ } y \in (x-\varepsilon,x+\varepsilon)\}
		\end{equation}
		be the value set of all local minima of $r$. The set $L$ is countable. Indeed, setting $I_{p,q} := \{ x \in (p,q) : r(x) \leq r(y) \text{ f.a.\ } y \in (p,q)\}$, we have $L = \bigcup _{(p,q) \in \mathbb{Q}^2} r(I_{p,q})$ where $r(I_{p,q})$ is either empty or a singleton. Since $u(X_\sigma,B_\sigma) = v(B_\sigma)$ and $X_\sigma = l \Rightarrow T_l \leq \sigma$ $\mathbb{P}$-a.s., we obtain
		\begin{align*}
		\xi[t > \sigma, X_{\sigma} \in L] &= \sum _{l \in L} \xi[t > \sigma, u(X_{\sigma},B_{\sigma}) = v(B_{\sigma}),  X_{\sigma} = l] \\
		&\leq \sum _{l \in L} \xi[t > \sigma \geq T_l, u(l,B_{\sigma}) = v(B_{\sigma})] 
		\end{align*} 
		and the r.h.s.\ is equal to $0$ by Lemma \ref{lemma:EquaToLeq}.
		
		It remains to show that $\xi[t > \sigma, X_{\sigma} \not \in L] = 0$. To this end, we define
		$[l]_n := \max \{ i/2^n : i \in \mathbb{N}, i/2^n \leq l \}$ for all $n \in \mathbb{N}$ and $l \geq 0$, and
		\begin{equation*}
		\sigma ^n := \inf \{t \geq 0: u([X_t]_n,B_t) = v(B_t) \}.
		\end{equation*}
		We claim that 
		\begin{equation} \label{eq:AuxClaim}
		\mathbb{P}\left[X_{\sigma} \not \in L, \sigma < \inf _{n \in \mathbb{N}} \sigma ^n\right] = 0.
		\end{equation}
		Admitting \eqref{eq:AuxClaim}, since for all $n \in \mathbb{N}$ the function $t \mapsto u([X_t]_n,B_t)$ is right-continuous, we have a.s.\ $u([X_{\sigma^n}]_n,B_{\sigma^n}) = v(B_{\sigma_n})$ and hence \eqref{eq:AuxClaim} yields
		\begin{align*}
		\xi[t > \sigma, X_{\sigma} \not \in L]  
		&\leq  \xi\left[ t > \inf _{n \in \mathbb{N}} \sigma^n \right] \\
		&\leq \sum _{n \in \mathbb{N}} \xi[ t > \sigma^n, u([X_{\sigma^n}]_n,B_{\sigma^n}) = v(B_{\sigma^n})] \\
		&= \sum _{n \in \mathbb{N}} \sum _{i = 0} ^{\infty} \xi \left[ t > \sigma^n, u([X_{\sigma^n}]_n,B_{\sigma^n}) = v(B_{\sigma^n}), \frac{i}{2^n} \leq X_{\sigma ^n} < \frac{i+1}{2^n} \right] \\
		&\leq \sum _{n \in \mathbb{N}} \sum _{i = 0} ^{\infty} \xi[ t > \sigma^n \geq i/2^n, u(i/2^n,B_{\sigma^n}) = v(B_{\sigma^n})].
		\end{align*}
		By Lemma \ref{lemma:EquaToLeq}, these summands are zero for all $n,i \in \mathbb{N}$.
		
		We are left with verifying \eqref{eq:AuxClaim}. By the definition of $L$ in \eqref{eq:DefOfL}, we see that for every pair $(l,x)$ where $l \not \in L$ and $x \in \mathbb{R}$ with $u(l,x) = v(x)$, there exists a sequence $(x_n)_{n \in \mathbb{N}}$ that converges to $x$ such that $u([l]_n,x_n) = v(x_n)$ for all $n \in \mathbb{N}$ large enough. Indeed, since $u(l,x) = v(x)$, it is $r(x) \leq l$ which leaves us with two cases: If $r(x) < l$, we just need to choose $n$ large enough such that $r(x) \leq [l]_n  \leq l$. If $r(x) = l \not \in L$, $x$ cannot be a local minimum of $r$, therefore there exists a sequence $(x_m)_{m \in \mathbb{N}}$ that converges to $x$ with $r(x_m) < l$ and we just need to choose an appropriate subsequence $(x_{m_n})_{n \in \mathbb{N}}$ such that $r(x_m) \leq [l]_{n_m} \leq l$.
		Thus, since  $u(X_\sigma,B_\sigma) = v(B_\sigma)$  $\mathbb{P}$-a.s., we obtain for $\mathbb{P}$-a.e.\ $\omega$
		\begin{align*}
		X_\sigma(\omega) \not \in L \quad &\Rightarrow \quad \forall \delta > 0 \, \exists n \in \mathbb{N} \, \exists y \in \mathcal{B}_\delta (B_{\sigma}(\omega)) \, : u([X_\sigma(\omega)]_n,y) = v(y) 
		\end{align*}
		where $\mathcal{B}_{\delta}(x)$ denotes the open ball of radius $\delta$ around $x$.
		Hence, for all $\varepsilon > 0$ we have
		\begin{equation} \label{eq:NastyIneq}
		\begin{split}
		&\mathbb{P}[\forall n \in \mathbb{N} \, \forall t \in (\sigma, \sigma + \varepsilon) : u([X_t]_n,B_t) < v(B_t), X_{\sigma} \not \in L ]   \\
		\leq \,&\mathbb{P}[\forall n \in \mathbb{N} \, \forall t \in (\sigma, \sigma + \varepsilon) : u([X_\sigma]_n,B_t) < v(B_t), X_{\sigma} \not \in L ] \\ 
		\leq \,& \mathbb{P}[\forall \delta  > 0 \, \exists y\in \mathcal{B}_\delta (B_{\sigma}) \, \forall t \in (\sigma, \sigma + \varepsilon) : B_t \neq y].
		\end{split}
		\end{equation}
		where we used the monotonicity of $u$ in the first component (cf.\ Lemma \ref{lemma:uCont}).
		By the strong Markov property and the continuity of Brownian motion, we can bound the last term in \eqref{eq:NastyIneq} by the sum of $\mathbb{P}[\forall t \leq \varepsilon : B_t \leq 0]$ and $\mathbb{P}[\forall t \leq \varepsilon : B_t \geq 0]$, and this is clearly $0$. Since $\varepsilon > 0$ is arbitrary, \eqref{eq:AuxClaim} is shown.
	\end{proof}
	
	Recall that $\hat{\tau} := \inf\{ t \geq 0 : u(X_t,B_t) = v(B_t)\}$ where $u(l,\cdot)$ is the potential function of $\Law(B_{\xi \land T_l})$ and $v$ is the potential function of $\nu$.
	
	\begin{corollary} \label{cor:Leq}
		We have a.s.\ $\xi[t \leq \hat{\tau}] = 1$. If $\xi$ is induced by an $\mathcal{F}$-stopping time $\tau$, we have $\tau \leq \hat{\tau}$.
	\end{corollary}
	
	\begin{proof}
		Since $u$ is continuous and $t \mapsto (X_t,B_t)$ is $\mathbb{P}$-a.s.\ right-continuous, we obtain  $\mathbb{P}[u(X_{\hat \tau}, B_{\hat \tau}) = v(B_{\hat \tau})] = 1$ and therefore we can apply Proposition \ref{prop:EquaToLeq}.
	\end{proof}

	\subsection{Proof of Theorem \ref{thm:MainEqui}} \label{ssec:ActualProof}
	
	Recall once again the properties of $(X_t)_{t \geq 0}$ and $(T_l)_{l \geq 0}$ formulated at the end of subsection \ref{ssec:adjoint}.
	Let $\xi$ be a RST which is a solution to $\SEP(\mu,\nu)$.
	Additionally to $\Law(B_{T_l \land \xi})$ (cf.\ \eqref{eq:DefRST2}), we introduce notation for the measures  \begin{equation} \label{eq:DefRST3}
	\begin{split}
	\Law(B_{\xi}; \xi \geq T_l) &:= ((\omega,t) \mapsto B_t(\omega))_{\#} \xi _{\vert \{t \geq T_l (\omega)\}} \quad \text{and} \\
	\Law(B_{T_l}; \xi \geq T_l) &:= ((\omega,t) \mapsto B_{T_l(\omega)}(\omega))_{\#} \xi _{\vert \{t \geq T_l (\omega)\}} 
	\end{split}
	\end{equation}
	
	The following Lemma \ref{lemma:ShadToSupp2} is the main observation that allows us to show in Lemma \ref{lemma:ShadToEqua} a counterpart to the upper bound stated in Corollary \ref{cor:Leq}. 
	
	\begin{lemma} \label{lemma:ShadToSupp2}
		Let $l \geq 0$ and suppose that $\xi$ satisfies 
		\begin{equation} \label{eq:ShaodwProp}
		\Law(B_{\xi}; \xi \geq T_l) = \shadow{\nu}{\Law(B_{T_l}; \xi \geq T_l)}.
		\end{equation}
		For all $x \in \mathbb{R}$,  if $u(l,x) < v(x)$,  $x \not \in \mathrm{supp}(\nu - \mathrm{Law}(B_{\xi}; \xi \geq T_l))$.
	\end{lemma}
	
	\begin{proof}
		Fix $l \geq 0$. By \eqref{eq:DefRST2} and \eqref{eq:DefRST3}, we have
		\begin{equation} \label{eq:DefinitionId}
		\Law(B_{T_l \land \xi}) - \Law(B_{T_l}; \xi \geq T_l) = \nu - \Law(B_{\xi}; \xi \geq T_l).
		\end{equation}
		Hence, Lemma \ref{lemma:PotfShad} and \eqref{eq:ShaodwProp} yield
		\begin{align*}
		v - u(l, \cdot) &= U_{\mathrm{Law}(B_{\xi};  \xi \geq T_l)} - U_{\mathrm{Law}(B_{T_l}; \xi \geq {T_l})} \\
		&= v - U_{\mathrm{Law}(B_{T_l}; \xi \geq T_l)} - \mathrm{conv} \left( v - U_{\mathrm{Law}(B_{T_l}; \xi \geq T_l)} \right).
		\end{align*}
		Let $x \in \mathbb{R}$ with $u(l,x) < v(x)$. By Lemma \ref{lemma:PropConv}, there exists an $\varepsilon > 0$ such that on the interval $[x- \varepsilon, x + \varepsilon]$ the function 
		\begin{align*}
		\mathrm{conv} \left( v - U _{\mathrm{Law}(B_{T_l}; \xi \geq T_l)} \right) 
		\end{align*}
		is affine. Rewriting with Lemma \ref{lemma:PotfShad} and \eqref{eq:DefinitionId}, we obtain that the function
		\begin{align*}
		u(l, \cdot) - U_{\mathrm{Law}(B_{T_l}; \xi \geq T_l)} = U_{\Law(B_{T_l \land \xi}) - \Law(B_{T_l}; \xi \geq T_l)} = U_{\nu - \Law(B_{\xi}; \xi \geq T_l)}
		\end{align*}
		is affine around $x$. Hence, Lemma \ref{lemma:PropPotf} yields that $x \not \in \mathrm{supp}(\nu - \Law(B_{\xi}; \xi \geq T_l))$.
	\end{proof}

	\begin{lemma} \label{lemma:ShadToEqua}
		If $\xi$ satisfies $\mathrm{Law}(B_{\xi}; \xi \geq T_l) = \shadow{\nu}{\mathrm{Law}(B_{T_l}; \xi \geq T_l)}$ for all $l \geq 0$, we have $\xi\left[ u(X_t, B_{t}) < v(B_{t}) \right] = 0$.
	\end{lemma}

	\begin{proof}	
		Let $\varepsilon > 0$. 
		Since both $u(0,\cdot)$ and $v$ are potential functions of probability measures with the same mass and barycenter, they are continuous and their difference vanishes at $\pm \infty$. Hence, there exists an $M_1 \in \mathbb{N}$ such that $v(x) - u(0,x) \leq \frac{\varepsilon}{2}$ for all $|x| \geq M_1$. On the compact interval $[-M_1,M_1]$ the monotone increasing sequence $(u(l,\cdot))_{l \geq 0}$ converges pointwise to $v$ (see Lemma \ref{lemma:uCont}). Dini's theorem yields that there exists $M_2 \in \mathbb{N}$ such that $\sup _{x \in [-M_1,M_1]} v(x) - u(l,x) \leq \frac{\varepsilon}{2}$ for all $l \geq M_2$.  Moreover, since $u$ is jointly continuous on the compact interval $[0,M_2] \times [-M_1,M_1]$,there exists an $n \in \mathbb{N}$ such that 
		\begin{equation*}
		\forall \, x \in \mathbb{R} \ \forall \, 0 \leq l \leq l' \leq l + \frac{1}{2^n} \, : \quad	u(l',x) - u(l,x) \leq \frac{\varepsilon}{2}.
		\end{equation*}
		For this $n$, we obtain 
		\begin{align*}
		\xi \left[ u(X_t, B_t) + \varepsilon \leq v(B_t) \right] &= \sum _{i = 1} ^{\infty} \xi \left[  u(X_t, B_t) + \varepsilon \leq v(B_t), \frac{i-1}{2^n} \leq X_t < \frac{i}{2^n} \right] \\
		&\leq \sum _{i = 1} ^{\infty} \xi \left[  u\left(\frac{i}{2^n}, B_t\right) < v(B_t),  X_t < \frac{i}{2^n} \right].
		\end{align*}
		For each $i \in \mathbb{N}$, the summands on the r.h.s.\ are $0$ because Lemma \ref{lemma:ShadToSupp2}  yields 
		\begin{align*}
		\xi \left[  u\left(\frac{i}{2^n}, B_t\right) < v(B_t),  X_t < \frac{i}{2^n} \right] =  \xi \left[  u\left(\frac{i}{2^n}, B_t \right) < v(B_t),  t < T_{\frac{i}{2^n}} \right] = 0
		\end{align*}
		Since $\varepsilon >0$ is arbitrary, the claim follows.
	\end{proof}

	\begin{lemma} \label{lemma:RBtoEqua}
		Let $\mathcal{R} \subset [0,\infty) \times \mathbb{R}$ be a closed barrier and $\tau$ an $\mathcal{F}$-stopping time. If $\tau = \inf \{t \geq 0: (X_t,B_t) \in \mathcal{R} \}$ $\mathbb{P}$-a.s., we have $\mathbb{P}[u(X_\tau,B_\tau) = v(B_{\tau})] = 1$.
	\end{lemma}
	
	\begin{proof}
		For all $(l,x) \in \mathcal{R}$, the Brownian motion $B$ cannot pass through $x$ on $[T_l \land \tau, \tau]$. Indeed, if $ t \in (T_l \land \tau, \tau]$ it is $X_t \geq l$ and since $\tau$ is by assumption the first time the process $(X_t,B_t)_{t \geq 0}$ hits the barrier $\mathcal{R}$, the Brownian motion is stopped at latest when it reaches $[l, \infty) \times \{x\} \subset \mathbb{R}$.
		Hence, we have  $(B_{\tau} - x)(B_{\tau \land T_l} - x) \geq 0$ $\mathbb{P}$-a.s., and thus we obtain $$u(l,x) = \mathbb{E}[|B_{T_l \land \tau }-x|] = \mathbb{E}[|B_{\tau} - x|] =  v(x).$$ 
		Since $u$ is continuous (cf.\  Lemma \ref{lemma:uCont}) and $t \mapsto (X_t,B_t)$ is right-continuous, we get $\mathbb{P}[u(X_{\tau},B_{\tau}) = v(B_{\tau})] = 1$.
	\end{proof}
	
	\begin{lemma} \label{lemma:InfimumToEqual}
		Let $\hat{\tau} := \inf\{t \geq 0 : u(X_t,B_t) = v(B_t)\}$. For all $l \geq 0$ we have $\mathbb{P}[\hat{\tau} < T_l, u(l,B_{T_l \land \hat{\tau}}) < v(B_{T_l \land \hat{\tau}})] = 0$.
	\end{lemma}
	
	\begin{proof}
		By Lemma \ref{lemma:uCont}, $u$ is continuous and  $t \mapsto (X_t,B_t)$ is $\mathbb{P}$-a.s.\ right-continuous, therefore we obtain from the definition of $\hat{\tau}$
		\begin{equation*}
		\mathbb{P}[u(X_{\hat \tau}, B_{\hat \tau}) = v(B_{\hat \tau})] = 1.
		\end{equation*}
		Let $l \geq 0$. Since $(X_t)_{t \geq 0}$ is adjoint to $(T_l)_{l \geq 0}$, we obtain
		\begin{equation*}
		\hat{\tau} < T_l \quad \Rightarrow \quad X_{\hat{\tau}} < l \quad \mathbb{P}\text{a.s.}
		\end{equation*}
		Since $u$ is also monotonously increasing in the first component and  it is  $B_{T_l \land \hat{\tau}} = B_{\hat{\tau}}$ on the set $\{\hat{\tau} < T_l\}$, the claim follows.
	\end{proof}

	Recall the definitions from the end of Subsection \ref{ssec:adjoint}:

	\begin{theorem} \label{thm:MainEqui}
		The following are equivalent:
		\begin{itemize}
			\item [(i)] There exists a closed barrier $\mathcal{R} \subset [0,\infty) \times \mathbb{R}$ such that $\xi$ is induced by the $\mathcal{F}$-stopping time $\tau := \inf \{ t \geq 0 : (X_t,B_t) \in \mathcal{R}\}$.
			\item [(ii)] For all $l \geq 0$ we have $\mathrm{Law}(B_{\xi}; \xi > T_l) = \shadow{\nu}{\mathrm{Law}(B_{\xi}; \xi > T_l)}$.
			\item [(iii)] $\xi$ is induced by the $\mathcal{F}$ stopping time $\hat{\tau} := \inf\{t \geq 0 : u(X_t,B_t) = v(B_t)\}$.
		\end{itemize}
	\end{theorem}
	
	\begin{proof}
		\textit{(i) $\Rightarrow$ (iii):} By Lemma \ref{lemma:RBtoEqua}, $\tau$ satisfies $\mathbb{P}[u(X_\tau,B_\tau) = v(B_{\tau})] = 1$ and thus we have $\mathbb{P}$-a.s. $\hat{\tau} \leq \tau$. The claim follows with Corollary \ref{cor:Leq}.
		
		\textit{(iii) $\Rightarrow$ (i):}  
		Lemma \ref{lemma:uCont} yields that $u$ is a jointly continuous function which is monotonously increasing in $l$. Hence, the set $\mathcal{R} := \{(l,x) \in [0, \infty) \times \mathbb{R} : u(l,x) = v(x)\}$ is a closed barrier.
		
		\textit{(ii) $\Rightarrow$ (iii):} 
		By Lemma \ref{lemma:ShadToEqua},
		$\xi[t \geq \hat \tau] \geq \xi[u(X_t,B_t) = v(B_t)] = 1$ and Corollary \ref{cor:Leq} yields that $\xi [t \leq \hat{\tau}] = 1$.
		
		\textit{(iii) $\Rightarrow$ (ii):} 
		Let $l \geq 0$.
		Since $\xi$ is induced by $\hat{\tau}$ and a solution to $\SEP(\mu,\nu)$, $\hat \tau - \hat{\tau} \land T_l$ is a solution of $\SEP(\Law(B_{\hat{\tau} \land T_l}),\nu)$ w.r.t.\ the Brownian motion $B'_s = B_{s + \hat{\tau} \land T_l}$.  Moreover, Lemma \ref{lemma:InfimumToEqual} yields 
		\begin{equation*}
		\hat{\tau} < T_l \Rightarrow u(l,B'_0) = v(B'_{0}) \quad \mathbb{P}\text{-a.s.}
		\end{equation*}
		Hence, by Corollary \ref{cor:ShadowOnEqualPart} it is
		\begin{align*}
		\Law(B_{\hat{\tau}}; \hat{\tau} \geq T_l) &= \Law(B'_{\hat{\tau} - \hat{\tau} \land T_l}; \hat \tau \geq T_l) \\
		&= \shadow{\nu}{\Law(B'_{0}; \hat \tau \geq T_l)} = \shadow{\nu}{\Law(B_{T_l}; \hat{\tau} \geq T_l)}. \qedhere
		\end{align*}
		\end{proof}

	\section{Proof of Proposition \ref{prop:ShiftedThm}}
	
	\begin{proof}
		Let $\tilde{\mathcal{F}}$ be the filtration defined by $\tilde{\mathcal{F}}_s := \mathcal{F}_{\sigma + s}$ and let $\tilde{B}$ be the process defined by $\tilde{B}_s := B_{\sigma + s}$. $\tilde{B}$ is an $\tilde{\mathcal{F}}$-Brownian motion. Moreover, $\tilde{\tau} := \tau - \sigma$ is an $\tilde{\mathcal{F}}$-stopping time because $\{\tilde{\tau} \leq s \} = \{\tau \leq \sigma + s\} \in \tilde{\mathcal{F}}_s$ for all $s \geq 0$. Clearly, we have $\tilde{B}_{\tilde{\tau}} = B_\tau$.
		
		Suppose (i) is satisfied. We set $\tilde{X}_s := X_{\sigma + s}$. Since $X$ is $\mathcal{F}$-adapted, $\tilde{X}$ is $\tilde{\mathcal{F}}$-adapted and, furthermore, we have
		\begin{equation*}
		\tilde{\tau} := \tau - \sigma = \inf \{ s \geq 0 : (\tilde{X}_s,\tilde{B}_s) \in \mathcal{R} \}.
		\end{equation*}
		Applying Theorem \ref{thm:intro} yields the existence of an $\tilde{\mathcal{F}}$-time-change $(\tilde{T}_l)_{l \geq 0}$ such that for all $l \geq 0$ we have
		\begin{equation*}
		\Law(\tilde{B}_{\tilde{\tau}}; \tilde{\tau} \geq \tilde{T}_l) = \shadow{\nu}{\Law(\tilde{B}_{\tilde{T}_l}; \tilde{\tau} \geq \tilde{T}_l}.
		\end{equation*}
		In particular, by Theorem \ref{thm:intro} we can choose $\tilde{T}_l = \inf \{ s \geq 0 : \tilde{X}_s \geq l\}$. Moreover, we set $T_l = \inf \{t \geq 0 : X_t \geq l\}$ and see that we have
		\begin{align*}
		\sigma + \tilde{T}_l &= \sigma + \inf \{ s \geq 0 : \tilde{X}_s \geq l\} = \sigma + \inf \{ s \geq 0 : X_{\sigma + s} \geq l\} \\
		&= \max\{\sigma, T_l \} 
		\end{align*}
		where the last equality follows from the fact that $X$ is monotonously increasing.
		We easily verify that a.s.
		\begin{equation} \label{eq:RelationTilde}
		\tilde{B}_{\tilde{T}_l} = B_{\sigma + \tilde{T_l}} = B_{\sigma \lor T_l} \quad \text{ and } \quad \{ \tilde{\tau} \geq \tilde{T}_l\} = \{\tau - \sigma \geq \tilde{T}_l \} \} = \{\tau \geq \sigma \lor T_l\}.
		\end{equation}
		Hence, for all $l \geq 0$ we obtain
		\begin{equation*}
		\Law(B_\tau; \tau \geq \sigma \lor T_l) = \shadow{\nu}{\Law(B_{\sigma \lor T_l}; \tau \geq \sigma \lor T_l )}.
		\end{equation*}

		Conversely, suppose that (ii) is satisfied. We set $\tilde{T}_l := \max\{0,T_l - \sigma\}$. Since $(T_l)_{l \geq 0}$ is an $\mathcal{F}$-time-change, $(\tilde{T}_l)_{l \geq 0}$ is an $\tilde{\mathcal{F}}$-time-change, and by definition we have $\sigma + \tilde{T}_l = \sigma \lor T_l$  such that \eqref{eq:RelationTilde} holds as well. In particular, we obtain
		\begin{align*}
		\Law(\tilde{B}_{\tilde{\tau}}; \tilde{\tau} \geq \tilde{T}_l) &= \Law(B_\tau; \tau \geq \sigma \lor T_l) \\
		&= \shadow{\nu}{\Law(B_{\sigma \lor T_l}; \tau \geq \sigma \lor T_l)} = \shadow{\nu}{\Law(\tilde{B}_{\tilde{T_l}}; \tilde{\tau} \geq \tilde{T}_l)}.
		\end{align*}
		Applying Theorem \ref{thm:intro} yields the existence of an $\tilde{\mathcal{F}}$-adapted stochastic process $(\tilde{X}_s)_{s \geq 0}$ and a closed barrier $\mathcal{R} \subset [0,\infty) \times \mathbb{R}$ such that 
		\begin{equation*}
		\tilde{\tau} = \inf \{ s \geq 0: (\tilde{X}_s,\tilde{B}_s) \in \mathcal{R}\}.
		\end{equation*} 
		In particular, by Theorem \ref{thm:intro} we can choose $\tilde{X}_s := \sup \{l \geq 0 : \tilde{T}_l \leq t\}$. Moreover, we set $X_t := \sup \{l \geq 0 : T_l \leq t\}$ and see that
		\begin{align*}
		\tilde{X}_s &= \sup \{l \geq 0: \tilde{T}_l \leq s\} = \sup \{l \geq 0: \max\{0, T_l - \sigma\} \leq s\} \\
		&= \sup \{l \geq 0 : T_l \leq \sigma + s\} = X_{\sigma + s}.
		\end{align*}
		Thus, we have a.s.
		\begin{align*}
		\inf \{ t \geq \sigma: (X_t, B_t) \in \mathcal{R} \} &= \sigma +  \inf\{ s \geq 0 : (X_{\sigma + s}, B_{\sigma + s}) \in \mathcal{R}\} \\
		&= \sigma + \tilde{\tau} = \tau. \qedhere
		\end{align*}
	\end{proof}

	\begin{remark} \label{rem:CondLM}
		For the left-monotone time-change $(T_l^{lm})_{l \geq 0}$, we have 
		\begin{equation*}
		\{\tau^i \geq \tau^{i-1} \lor T_l ^{lm} \} = \{\tau^i \geq \tau^{i-1}, T^{lm}_l = 0\} = \{ B_0 \leq q_l\}
		\end{equation*}
		where $q_l :=- \ln (l)$ and $T_l ^{lm} = 0$ on this set. Thus, the stopping times $\tau ^1 \leq ... \leq \tau ^n$ are $(T_l ^{lm})_{l \geq 0}$-shadow-residual if and only if
		\begin{align*}
		\Law(B_{\tau ^i}; B_0 \leq q_l) &= \Law(B_{\tau ^i}; \tau ^i \geq \tau^{i-1}  \lor T_l ^{lm})  \\
		&= \shadow{\nu _i}{\Law(B_{\tau ^{i-1} \lor T_l ^{lm}}; \tau ^i \geq \tau^{i-1}  \lor T_l ^{lm})} \\
		&= \shadow{\nu _i}{\Law(B_{\tau ^{i-1}}; B_0 \leq q_l)}
		\end{align*}
		for all $1 \leq i \leq n$. Applying this inductively, these stopping times are shadow-residual if and only if
		\begin{equation*}
		\Law(B_{\tau ^i}; B_0 \leq q_l) = \shadow{\nu _i}{ ... \, \shadow{\nu_1}{\Law(B_{0}; B_0 \leq q_l)}} =: \shadow{\nu_1,...,\nu _i}{\Law(B_{0}; B_0 \leq q_l)}
		\end{equation*}
		for all $1 \leq i \leq n$. This is the obstructed shadow defined by Nutz-Stebegg-Tan in \cite{NuStTa17}. Hence, $(\tau ^1, ... , \tau^n)$ is the multi-marginal lm-solution if and only if the joint distribution $(B_0,B_{\tau ^1}, ... , B_{\tau ^n})$ is the mutliperiod left-monotone transport.
	\end{remark}

	\section{Proof of Proposition \ref{prop:Interpolation}} 
	\label{sec:Interpolation}
	
	In this subsection we suppose that $\Omega = C([0,\infty))$ is the path space of continuous functions and that $\mathbb{P}$ is a probability measure on the path space such that the canonical process $B: \omega \mapsto \omega$ is a Brownian motion with $\Law_{\mathbb{P}}(B_0) = \mu$. Moreover, we denote by $\theta$ the shift operator on $\Omega$, i.e.\ $\theta_r : (\omega_s)_{s \geq 0} \mapsto (\omega_s)_{s \geq r}$ for all $r \geq 0$.
	
	\subsection{Concatenation Method} 
	To simplify notation, we say that a finite stopping time $\tau$ is shadow-residual w.r.t.\ a time-change $(T_l)_{l \geq 0}$ if for all $l \geq 0$ we have
	\begin{equation*}
	\Law(B_\tau; \tau \geq T_l) = \shadow{\Law(B_\tau)}{\Law(B_{T_l}; \tau \geq T_l)}.
	\end{equation*}
	This is precisely the condition in part (ii) of Theorem \ref{thm:intro}.

	\begin{lemma} \label{lemma:CombinedStoppingTime}
		Let $\tau$ and $\sigma$ be two $\mathcal{F}$-stopping times such that $\tau$ is finite. The random variable $\tau + \sigma \circ \theta_{\tau}$ is a again a $\mathcal{F}$ stopping time.
		\end{lemma}
	
	\begin{proof}
		If $\tau$ takes only values in the countable set $A \subset [0, \infty)$, for all $s \geq 0$ we obtain
		\begin{equation*}
		\{\tau + \sigma \circ \ \theta_\tau \leq s\} = \bigcup _{k \in A \cap [0,t]} \{\sigma \circ \theta _k \leq t - k\} \in \mathcal{F}_{k + (t-k)} = \mathcal{F}_t.
		\end{equation*}
		A general $\tau$ can be approximated by discrete stopping times. 
	\end{proof}

	\begin{corollary} \label{lemma:NestingPrep}
		Let $(T_l)_{l \geq 0}$ be a finite $\mathcal{F}$-time-change, $(S_l)_{l \geq 0}$ a $\mathcal{F}$-time-change and $\lambda > 0$.
		The family $(R_l)_{l \geq 0}$ defined by 
		\begin{equation*}
		R_l := T_{l \land \lambda} + (S_{l-\lambda} \circ \theta _{T_\lambda}) \1 _{\{l \geq \lambda\}} = \begin{cases}
		T_l & l < \lambda \\
		T_\lambda + S_{l - \lambda} \circ \theta _{T_\lambda} &l \geq \lambda
		\end{cases}
		\end{equation*}
		is an $\mathcal{F}$-time-change. 
		If additionally, both $(T_l)_{l \geq 0}$ and $(S_l)_{l \geq 0}$ are left-continuous, $T_0 = S_0 = 0$, $T_\infty = S_\infty = + \infty$ and 
		$\mathbb{P}[\lim _{k \downarrow l} T_k = T_l] =\mathbb{P}[\lim _{k \downarrow l} S_k = S_l] = 1$  for all $l \geq 0$,
		$(R_l)_{l \geq 0}$ satisfies these four properties as well.
	\end{corollary}
	
	\begin{lemma} \label{lemma:Nesting}
		Suppose we are in the setting of Lemma \ref{lemma:NestingPrep}.
		Additionally assume that $\tau$ is a solution of $\SEP(\mu,\nu)$ which is shadow residual w.r.t.\ $(T_l)_{l \geq 0}$. If $\sigma$ is a $\mathcal{F}$-stopping time such that $\sigma$ is a solution to $\SEP(\mathrm{Law}(B_{\tau \land T_\lambda}),\nu)$, then 
		\begin{equation*}
		\rho := \tau \land T_{\lambda} + \sigma \circ \theta _{T_\lambda \land \tau} 
		\end{equation*}
		is a $\mathcal{F}$-stopping time and a solution to $\SEP(\mu,\nu)$ which is shadow residual w.r.t.\ $(R _l)_{l \geq 0}$.
	\end{lemma}

	\begin{proof}
		
		We set $\tilde{\mathcal{F}}_s = \mathcal{F}_{s + \tau \land T_\lambda}$, $\tilde{B} := B \circ \theta _{\tau \land T_\lambda}$, $\tilde{\sigma} := \sigma \circ \theta_{\tau \land T_\lambda}$ and $\tilde{S}_l := S_l \circ \theta_{T_\lambda \land \tau}$. $\tilde{\sigma}$ is a stopping time w.r.t.\ the filtration generated by $\tilde{B}$. We also have $\Law(\tilde{B}_{\tilde{\sigma}}) = \nu$ and $\tilde{\sigma}$ is $(\tilde{S}_l)_{l \geq 0}$-shadow-residual.
		
		\texttt{STEP 1:}
		We have $\Law(B_{\rho}) = \Law(B_{\tau \land T_\lambda + \tilde{\sigma}}) =  \mathrm{Law}(\tilde{B}_{\tilde{\sigma}}) = \nu$. By Lemma \ref{lemma:CombinedStoppingTime}, $\rho$ is a $\mathcal{F}$-stopping time and $(B_{s \land \rho})_{s \geq 0}$ is uniformly integrable because $\tau$ and $\sigma$ are solutions to $\SEP(\mu,\nu)$ and $\SEP(\Law(\tilde{B}_0), \nu)$. Thus, $\rho$ is a solution to $\SEP(\mu,\nu)$. Moreover, we claim that we can represent $\rho$ as
		\begin{equation} 	\label{eq:Step1}
		\rho = \begin{cases}
		\tau & \tau < T_\lambda \\
		T_\lambda + \tilde{\sigma} & \tau \geq T_\lambda
		\end{cases} \quad \mathbb{P}\text{-a.s.}.
		\end{equation}
		Indeed, since  $\tau$ is $(T_l)_{l \geq 0}$-shadow-residual, by Theorem  \ref{thm:MainEqui} (iii) and Lemma \ref{lemma:InfimumToEqual}, we have  $\tau < T_\lambda \Rightarrow u(\lambda,\tilde{B}_0) = v(\tilde{B}_0)$ $\mathbb{P}$-a.s.\ where $u(l,\cdot) := U_{\Law(B_{T_l \land \tau})}$ and $v = U_{\nu}$ for all $l \geq 0$. Thus, we get
		\begin{align*}
		\mathbb{P}[\tilde{\sigma} > 0, \tau < T_\lambda] &\leq \mathbb{P}[\tilde{\sigma} > 0, u(\lambda, \tilde{B}_0) = v(\tilde{B}_0)] \\
		&= \mathbb{P}[\tilde{\sigma} > 0, U_{\Law(\tilde{B}_0)}(\tilde{B}_0) 
		= U_{\nu}(\tilde{B}_0)].
		\end{align*}
		and the r.h.s. is equal to $0$  because $\tilde{\sigma}$ is a $\tilde{\mathcal{F}}$-stopping-time that solves $\SEP(\Law(\tilde{B}_0),\nu)$ (cf.\ Lemma \ref{lemma:EquaToLeq}).
		It remains to show that $\rho$ is $(R_l)_{l \geq 0}$-shadow-residual. We split this up in the cases $l \geq \lambda$ and $l < \lambda$.
		
		\texttt{STEP 2:}
		Suppose $l \geq \lambda$. Since by \texttt{STEP 1}  $\{ \rho \geq R_l \} = \{\tilde{\sigma} \geq \tilde{S}_{l - \lambda},\tau \geq T_\lambda\}$ $\mathbb{P}$-a.s.\ and $\tilde{\sigma}$ is $(\tilde{S}_l)_{l \geq 0}$ shadow-residual, Lemma \ref{lemma:ShadowAssz} yields
		\begin{align*}
		&\Law(B_\rho; \rho \geq R_l) + \Law(\tilde{B}_{\tilde{\sigma}}; \tilde{\sigma} \geq \tilde{S}_{l-\lambda}, \tau < T_\lambda) \\ &\quad =  \Law(\tilde{B}_{\tilde{\sigma}}; \tilde{\sigma} \geq \tilde{S}_{l - \lambda})
		= \shadow{\nu}{\Law(\tilde{B}_{\tilde{S}_{l - \lambda}}; \tilde{\sigma} \geq \tilde{S}_{l - \lambda} )}
		\\
		& \quad = \shadow{\nu}{\Law(B_{R_l}; \rho \geq R_l)} \\  
		& \hspace{2cm} + \shadow{\nu - \shadow{\nu}{\Law(\tilde{B}_{\tilde{S}_{l-\lambda}}; \tilde{\sigma} \geq \tilde{S}_{l-\lambda}, \tau \geq T_\lambda)}}{\Law(\tilde{B}_{\tilde{S}_{l-\lambda}}; \tilde{\sigma} \geq \tilde{S}_{l-\lambda}, \tau < T_\lambda)}
		\end{align*}
		Thus, we obtain $\Law(B_\rho; \rho \geq R_l) = \shadow{\nu}{\Law(B_{R_l}; \rho \geq R_l)}$ if we show
		\begin{align}
		&\Law(\tilde{B}_{\tilde{\sigma}}; \tilde{\sigma} \geq \tilde{S}_{l-\lambda}, \tau < T_\lambda) =  \Law(\tilde{B}_{\tilde{S}_{l-\lambda}}; \tilde{\sigma} \geq \tilde{S}_{l-\lambda}, \tau < T_\lambda) \quad \text{and} \label{eq:Step2Eq1}\\
		&\Law(\tilde{B}_{\tilde{S}_{l-\lambda}}; \tilde{\sigma} \geq \tilde{S}_{l-\lambda}, \tau < T_\lambda) \leq_+ \nu - \shadow{\nu}{\Law(\tilde{B}_{\tilde{S}_{l-\lambda}}; \tilde{\sigma} \geq \tilde{S}_{l-\lambda}, \tau \geq T_\lambda)}. \label{eq:Step2Eq2}
		\end{align}
		By \texttt{STEP} 1 it is $\tilde{\sigma} = 0$ on $\{\tau < T_\lambda\}$ and therefore \eqref{eq:Step2Eq1} follows immediately. Moreover, Lemma \ref{lemma:ShadowAssz} yields
		\begin{align*}
		\shadow{\nu}{\Law(\tilde{B}_{\tilde{S}_{l-\lambda}}; \tilde{\sigma} \geq \tilde{S}_{l-\lambda}, \tau \geq T_\lambda)} &\leqp 	\shadow{\nu}{\Law(\tilde{B}_{\tilde{S}_{l-\lambda} \land \tilde{\sigma}}; \tau \geq T_\lambda)}
		\end{align*}
		On the one hand, by the definition of the shadow we have
		\begin{equation} \label{eq:Aux2}
		\begin{split}
		\shadow{\nu}{\Law(\tilde{B}_{0}; \tau \geq T_\lambda)} &\leqc  \shadow{\nu}{\Law(\tilde{B}_{\tilde{S}_{l-\lambda} \land \tilde{\sigma}}; \tau \geq T_\lambda)} \\
		&\leq_c \shadow{\nu}{\Law(\tilde{B}_{\tilde{\sigma}}; \tau \geq T_\lambda)} = \Law(\tilde{B}_{\tilde{\sigma}}; \tau \geq T_\lambda)
		\end{split}
		\end{equation}
		because $0 \leq \tilde{S}_{l-\lambda} \land \tilde{\sigma} \leq  \tilde{\sigma}$ and $\Law(\tilde{B}_{\tilde{\sigma}}; \tau \geq T_\lambda) \leqp \nu$. On the other hand, since $u(\lambda, \tilde{B}_0) = v(\tilde{B}_0)$ on $\{\tau < T_\lambda\}$ by \texttt{STEP 1}, Corollary \ref{cor:ShadowOnEqualPart} yields
		\begin{equation*}
		\shadow{\nu}{\Law(\tilde{B}_{0}; \tau \geq T_\lambda)} = \Law(\tilde{B}_{\tilde{\sigma}}; \tau \geq T_\lambda).
		\end{equation*}
		Thus, we have equality in \eqref{eq:Aux2} which implies 
		\begin{equation*}
		\shadow{\nu}{\Law(\tilde{B}_{\tilde{S}_{l-\lambda}}; \tilde{\sigma} \geq \tilde{S}_{l-\lambda}, \tau \geq T_\lambda)} \leqp \nu - \Law(\tilde{B}_{\tilde{\sigma}}; \tau < T_\lambda) 
		\end{equation*}
		and thereby \eqref{eq:Step2Eq2}.

		\texttt{STEP 3:}
		Now suppose $l < \lambda$.  Since $R_\lambda = T_\lambda$ and $\{\rho \geq R_\lambda\} = \{\tau \geq T_\lambda\}$ $\mathbb{P}$-a.s., by \texttt{STEP 2} we have
		\begin{align*}
		\Law(B_\rho; \rho \geq R_\lambda) &= \shadow{\nu}{\Law(B_{R_\lambda}; \rho \geq R_\lambda)}  = \shadow{\nu}{\Law(B_{T_\lambda}; \tau \geq T_\lambda)} \\
		&= \Law(B_\tau; \tau \geq T_\lambda)
		\end{align*}
		because $\tau$ is $(T_l)_{l \geq 0}$ shadow residual.
		In particular, we get
		\begin{align*}
		\Law(B_\rho; \rho \geq R_l) &= \Law(B_\rho; \rho \geq R_l, \tau < T_\lambda ) + \Law(B_\rho; \rho \geq R_l, \tau \geq T_\lambda ) \\
		&= \Law(B_\tau; T_\lambda > \tau \geq T_l) + \Law(B_\rho; \rho \geq T_\lambda) \\
		&= \Law(B_\tau; \tau \geq T_l) \\
		&= \shadow{\nu}{\Law(B_{T_l}; \tau \geq T_l)} \\
		&= \shadow{\nu}{\Law(B_{R_l}; \rho \geq R_l)}
		\end{align*}
		because $\{\rho \geq R_l \} = \{\tau \geq T_\lambda\} \cup \{T_\lambda > \tau \geq T_\lambda\} = \{\tau \geq T_l\}$.
	\end{proof}

	\subsection{Robustness of LM-Embedding}

	\begin{lemma} \label{lemma:SEPcompact}
		Let $(\mathbb{P}_n)_{n \in \mathbb{N}}$ be a sequence of probability measures on $\Omega$ such that $B$ is a Brownian motion under $\mathbb{P}_n$, $(\nu_n)_{n \in \mathbb{N}}$ a sequence of probability measures on $\mathbb{R}$ and $(\xi^n)_{n \in \mathbb{N}}$ a sequence of RST w.r.t.\ $\mathbb{P}_n$ that are a solution to $\SEP(\Law_{\mathbb{P}_n}(B_0), \nu _n)$ for all $n \in \mathbb{N}$. If $(\mathbb{P}_n)$ converges weakly to $\mathbb{P}$ and $(\nu_n)_{n \in \mathbb{N}}$ converges to the probability measure $\nu$ under $\TO$, there exists a weakly convergent subsequence of $(\xi ^n)_{n \in \mathbb{N}}$. Moreover, the limit of every convergent subsequence is a RST w.r.t.\ $\mathbb{P}$ that solves $\SEP(\Law_\mathbb{P}(B_0),\nu)$.
	\end{lemma}
	
	\begin{proof}
		Let $\varepsilon > 0$.
		Since $(\mathbb{P}_n)_{n \in \mathbb{N}}$ converges weakly, there exists a compact set $K_{\varepsilon} \subset \Omega$ such that $\mathbb{P}^n[K_\varepsilon] > 1 - \varepsilon$ for all $n \in \mathbb{N}$. 
		Since $(\nu_n)_{n \in \mathbb{N}}$ converges in $\TO$, by Lemma \ref{lemma:T1Conv} there exists $\eta \in \MO$ such that $\int \varphi \de \nu _n \leq \int \varphi \de \eta$ for all $n \in \mathbb{N}$ and non-negative convex functions $\varphi$. Moreover, by the Theorem of de la Vallee-Poussin there exists a non-negative convex function $V \in C^2(\mathbb{R})$ with $V'' \geq C > 0$ such that $\int_\mathbb{R} V \de \eta < \infty$. For all $s \geq 0$ and $n \in \mathbb{N}$ we have
		\begin{equation*}
		\xi ^n [t \geq s] = \overline{\mathbb{P}}[\overline{\tau}^{\xi^n} \geq s] \leq \frac{\overline{\mathbb{E}}[\overline{\tau} ^{\xi^n}]}{s} \leq \frac{\overline{\mathbb{E}}[V(\overline{B}_{\overline{\tau} ^{\xi^n}})]}{s} \leq \frac{1}{Cs} \int_{\mathbb{R}} V \de \eta
		\end{equation*}
		where we used the notation of Lemma \ref{lemma:ReprRST}, the Markov inequality and Ito's formula.
		Hence, there exists $s_\varepsilon > 0$ such that $\xi^n[t \leq s_\varepsilon] > 1 - \varepsilon$ for all $n \in \mathbb{N}$.
		Then the mass of the compact set $K_\varepsilon \times [0,s_\varepsilon]$ under $\xi ^n$ is strictly greater than $1- 2\varepsilon$ for all $n \in \mathbb{N}$. Hence,  the set $\{\xi ^n : n \in \mathbb{N} \}$ is tight. 
		By Prokhorovs Theorem there exists a weakly convergent subsequence. We denote the limit by $\xi$.
		
		Since the set of RST is  closed under weak convergence (cf.\ \cite[Corollary 3.10]{BeCoHu17}), $\xi$ is a RST stopping time w.r.t.\ $\mathbb{P}$. Moreover, $(\omega,t) \mapsto \varphi(\omega_t)$ is a continuous and bounded function on $\Omega \times [0, \infty)$  for all $\varphi \in C_b(\mathbb{R})$, and therefore $\Law(B_\xi) = \nu$. 
		It remains to show that $(B_{\xi \land t})_{t \geq 0}$ is uniformly integrable. Since $|x| \1 _{|x| \geq K} \leq |x- K/2| + |x + K/2| - K$, we get
		\begin{equation*}
		\mathbb{E} \left[ |B_{\xi \land s}| \1 _{\{|B_{\xi \land s}| \geq K\}} \right]  \leq U_{\Law(B_{\xi \land s})}(K/2) + U_{\Law(B_{\xi \land s})} (-K/2) - K
		\end{equation*}
		for all $s,K \geq 0$.
		Moreover, since $\xi ^n$ converges weakly to $\xi$ and $g_m$ defined by $g_m(y) := \min\{|y|,m\}$ is a continuous and bounded function, we obtain for all $x \in \mathbb{R}$ with monotone and dominated convergence
		\begin{align*}
		U_{\Law(B_{\xi \land t})}(x) &= \sup_{m \in \mathbb{N}} \lim_{n \rightarrow \infty} \int _{\Omega \times [0, \infty)} g_m(\omega_{t \land s} - x) \de \xi ^n(\omega,t) \\
		&\leq \sup _{n \in \mathbb{N}} \mathbb{E}\left[ |B_{\xi^{n} \land s}| \right] \leq \sup _{n \in \mathbb{N}} \mathbb{E}\left[ |B_{\xi^{n}}| \right]  = \sup _{n \in \mathbb{N}} U_{\nu_n}(x) \leq U_{\eta}(x)
		\end{align*}
		where we used that $\Law(B_{\xi ^n \land t})_{t \geq 0}$ is uniformly integrable for all $n \in \mathbb{N}$ .
		Thus, using the asymptotic behaviour of potential functions, it is
		\begin{equation*}
		\lim _{K \rightarrow \infty} \sup _{t \geq 0} \mathbb{E} \left[ |B_{\xi \land t}| \1 _{\{|B_{\xi \land t}| \geq K\}} \right] \leq \limsup _{K \rightarrow \infty}  U_{\eta}\left(-K/2\right) + U_{\eta}(K/2) - K = 0.
		\end{equation*}
		and the claim follows.
	\end{proof}
	
	\begin{lemma} \label{lemma:StabilityLM}
		Let $(\nu_n)_{n \in \mathbb{N}}$ be a sequence of probability measures on $\mathbb{R}$, $(\mathbb{P}^n)_{n \in \mathbb{N}}$ be a sequence of probability measures on $\Omega$ such that $B$ is a Brownian motion with initial distribution $\mu_n$  under $\mathbb{P}^n$ and $(\xi^n)_{n \in \mathbb{N}}$ a sequence of corresponding RST which are lm-monotone solutions to $\SEP(\mu_n, \nu_n)$. 
		If $(\mathbb{P}_n)_{n \in \mathbb{N}}$ converges weakly to $\mathbb{P}$ whose initial distribution $\mu$ is atomless and $\nu_n$ converges to $\nu$ in $\TO$, the sequence $(\xi ^ n)_{n \in \mathbb{N}}$ converges weakly to a RST $\xi$ w.r.t.\ $\mathbb{P}$ which is a lm-monotone solution to $\SEP(\mu,\nu)$
	\end{lemma}
	
	\begin{proof}
		By Lemma \ref{lemma:SEPcompact}, any subsequence of $(\xi ^n)_{n \in \mathbb{N}}$ has itself a convergent subsequence and the limit is a solution to $\SEP(\mu,\nu)$. If we show that this limit is $(T_l^{lm})_{l \geq 0}$-shadow-residual, by uniqueness 
		(see \cite[Lemma 4.3]{BeHeTo17}), it has to be the unique left-monotone solution to $\SEP(\mu,\nu)$ and the claim follows.
		
		For simplicity, we denote the convergent subsequence of a given subsequence again by $(\xi ^n)_{n \in \mathbb{N}}$ and the limit by $\xi$. Since $\xi^n$ is $(T_l ^{lm})$-shadow-residual, the probability measure $\Law(B_0,B_{\xi ^n})$ is the left-curtain coupling of $\mu_n$ and $\nu_n$. Indeed, as in Remark \ref{rem:CondLM} we have for all $l \geq 0$
		\begin{equation} \label{eq:LmStability} 
		\begin{split}
		\Law(B_{\xi^n}; B_0 \leq - \ln(l)) &= \Law(B_{\xi ^n}; \xi ^n \geq T_l ^{lm}) \\ &= \shadow{\nu ^n}{\Law(B_0; \xi ^n \geq T_l ^n)} 
		= \shadow{\nu _n}{\Law(B_0; B_0 \leq - \ln(l))}
		\end{split}
		\end{equation}
		because $\{t \geq T_l ^{lm}\} = \{T_l ^{lm} = 0\} = \{B_0 \leq - \ln(l) \}$ (where $-\ln(0) := - \infty$). As shown in \cite[Theorem 2.16]{Ju14} (and also as a consequence of stability of martingale optimal transport \cite[Theorem 1.1]{BaPa19})), the left-curtain coupling is stable under weak convergence, i.e.\ the weak limit $\Law(B_0,B_\xi)$ of $(\Law(B_0,B_{\xi^n}))_{n \in \mathbb{N}}$ is the left-curtain coupling of $\mu$ and $\nu$. Thus, analogous to \eqref{eq:LmStability}, $\xi$ is $(T_l^{lm})_{l \geq 0}$-shadow-residual.
	\end{proof}

	\subsection{Application}
	
	Fix $\mu \leqc \nu$, let $\tau ^r$ be the Root solution to $\SEP(\mu,\nu)$ and $(T_l^r)_{l \geq 0}$ the Root time-change. 
	Let $\lambda > 0$. We set $ \tilde B^{\lambda} := B \circ \theta _{T^{r} _{\lambda} \land \tau ^{r}}$. By the strong Markov property, $ \tilde B^{\lambda}$ is a Brownian motion and there exists a left-monotone solution $\sigma ^{\lambda}$ of $\SEP(\mathrm{Law}(\tilde B^{\lambda}_0),\nu)$. We define
	\begin{align*}
	\tau ^{\lambda} :=& \tau ^{r} \1 _{\{ \tau ^{r} < T^{r}_{\lambda} \}} + \sigma^\lambda \circ \theta_{(T^{r} _{\lambda} \land \tau ^{r})} \1 _{\{ \tau ^{r} \geq  T^{r} _{\lambda} \}} \\
	=& \tau^r \1_{\{\tau ^r < \lambda \}} + \sigma^\lambda \circ \theta _{\lambda} \1 _{\{\tau ^r \geq \lambda\}} .
	\end{align*}
	By Lemma \ref{lemma:Nesting}, $\tau ^\lambda$ is a solution to $\SEP(\mu,\nu)$ which is shadow residual w.r.t.\ the time-change $(T_l ^\lambda)_{l \geq 0}$ defined as
	\begin{equation*}
	T^{\lambda} _l :=
	T^{r}_{l \land \lambda} + (T^{lm} _{l - \lambda} \circ \theta _{T^{r}_{\lambda}}) \1_{\{l \geq \lambda\}} = \begin{cases}
	l &  l < \lambda \\
	l & \exp(-B_\lambda) + \lambda \geq l \geq \lambda \\
	+ \infty  & \exp(-B_\lambda) + \lambda < l, l \geq \lambda
	\end{cases}.
	\end{equation*}
	Thus, by Theorem \ref{thm:MainEqui}, there exists a barrier $\mathcal{R}^\lambda$ such that 
	\begin{equation*}
	\tau ^{\lambda} = \inf \{t \geq 0 : (X_t ^\lambda, B_t) \in \mathcal{R}^\lambda\}
	\end{equation*}
	where $X^{\lambda}$ is defined as
	\begin{align*}	
	X_t ^\lambda := \sup \{l \geq 0 : T_l^\lambda \leq t\} = \begin{cases}	
	t & t < \lambda \\
	\lambda + \exp(-B_0)   & t \geq \lambda	\end{cases}.
	\end{align*}
	To complete the proof of Proposition \ref{prop:Interpolation}, it remains to show the  convergence of $\tau ^\lambda$ to $\tau ^r$ and $\tau ^{lm}$ as randomized stopping times as $\lambda$ tends to $+\infty$ and $0$. This is covered by Lemma \ref{lemma:ConvToRoot} and Lemma \ref{lemma:ConvToLM}.

	\begin{lemma} \label{lemma:ConvToRoot}
		The sequence $(\tau _\lambda)_{\lambda > 0}$ converges a.s.\ to $\tau ^r$ as $\lambda$ tends to $+ \infty$. In particular, $\Law(B,\tau ^\lambda)$ converges weakly to $\Law(B,\tau ^r)$.
	\end{lemma}
	
	\begin{proof}
		Since $\tau^r < + \infty$ and $T_\lambda^r \rightarrow \infty$, we have
		\begin{equation*}
		\lim _{\lambda \rightarrow + \infty} \tau ^\lambda = \lim _{\lambda \rightarrow + \infty} \left( \tau ^{r} \1 _{\{ \tau ^{r} < T^{r}_{\lambda} \}} + \sigma^\lambda \circ \theta_{(T^{r} _{\lambda} \land \tau ^{r})} \1 _{\{ \tau ^{r} \geq  T^{r} _{\lambda} \}} \right) = \tau ^{r} \quad \mathbb{P}-\text{a.s.}
		\end{equation*}
		The weak convergence of $\Law(B,\tau ^\lambda)$ to $\Law(B,\tau ^r)$ follows immediately.
	\end{proof}

	\begin{lemma} \label{lemma:CompSupp}
		For all $\varphi \in C_c(\Omega \times [0, \infty))$ we have
		\begin{equation*}
		\lim_{\lambda \rightarrow 0} \mathbb{E}\left[ \left\vert \varphi(B, \tau ^\lambda) - \varphi (\tilde{B}^\lambda, \tilde{\sigma}^\lambda) \right\vert \1 _{\{\tau^r > 0\}}\right] = 0
		\end{equation*}
		where $\tilde{B}^\lambda := B \circ \theta_{\tau ^r \land \lambda}$ and $\tilde{\sigma}^\lambda := \sigma^\lambda \circ \theta_{\tau ^r \land \lambda}$ for all $\lambda > 0$.
	\end{lemma}
	
	\begin{proof}
		For all $\lambda > 0$ we define the map $\Theta ^\lambda$ on $\Omega \times [0,\infty)$ by 
		\begin{equation*}
		\Theta ^\lambda : (\omega,t) \mapsto (\omega \circ \theta_\lambda, \max\{t-\lambda,0\}).
		\end{equation*}
		A compatible metric on the Polish space $\Omega \times [0, \infty)$ is given by
		\begin{equation*}
		d((\omega,t),(\omega',t')) := |t-t'| + \sum _{n \in \mathbb{N}}  2^{-n}\sup _{s \in [0,n]} |\omega_s - \omega'_s| 
		\end{equation*}
		and under this metric $\Theta ^\lambda$ is $2$-Lipschitz continuous for all $\lambda \in (0,1)$. Moreover, since $\lim _{\lambda \rightarrow 0}\Theta ^\lambda(\omega,t) = (\omega,t)$ for all  $(\omega,t) \in \Omega \times [0, \infty)$, $\Theta^\lambda$ converges uniformly on compact sets to the idenity on $\Omega \times [0, \infty)$. Thus, we have
		\begin{align*}
		&\lim _{\lambda \rightarrow 0} \mathbb{E}\left[ \left\vert \varphi(B,\lambda + \tilde{\sigma} ^\lambda) - \varphi(\tilde{B}^\lambda,\tilde{\sigma} ^\lambda)  \right\vert \right] \\
		& \quad = \lim _{\lambda \rightarrow 0} \mathbb{E}\left[ \varphi(B,\lambda + \tilde{\sigma} ^\lambda) - \varphi(\Theta^\lambda(B,\lambda + \tilde{\sigma} ^\lambda)) \right] = 0.
		\end{align*} 
		By substituting the definition of $\tau^\lambda$, we obtain the estimate
		\begin{align*}
		&\mathbb{E}\left[ \left\vert \varphi(B, \tau ^\lambda) - \varphi (\tilde{B}^\lambda, \tilde{\sigma}^\lambda) \right\vert \1 _{\{\tau^r > 0\}}\right] \\
		&\quad \leq 2 ||\varphi||_\infty \mathbb{P}[0 <\tau ^r < \lambda] + \mathbb{E}\left[ \left\vert \varphi(B,\lambda + \tilde{\sigma} ^\lambda) - \varphi(\tilde{B}^\lambda,\tilde{\sigma} ^\lambda) \right\vert \right]
		\end{align*}
		and therefore the claim follows.
	\end{proof}

	\begin{lemma} \label{lemma:ConvToLM}
		The sequence $(\Law(B,\tau ^\lambda))_{\lambda > 0}$ converges  weakly to $\Law(B,\tau ^{lm})$ as $\lambda$ tends to $0$.
	\end{lemma}
	
	\begin{proof}
		On the set $\{\tau ^r = 0\}$, $U_\mu(0) = u^r(0,B_0) = v(B_0) = U_\nu(B_0)$ and thus $\tau ^{lm} = 0 = \tau ^r$.
		Hence, in conjunction with Lemma \ref{lemma:CompSupp} we obtain for all $\varphi \in C_c(\Omega \times [0, \infty))$
		\begin{equation} \label{eq:ConvPHI}
		\lim_{\lambda \rightarrow 0} \mathbb{E}\left[ \left\vert \varphi(B, \tau ^\lambda) - \varphi (\tilde{B}^\lambda, \tilde{\sigma}^\lambda) \right\vert \right] = 0
		\end{equation}
		where $\tilde{B}^\lambda := B \circ \theta_{\tau ^r \land \lambda}$ and $\tilde{\sigma}^\lambda := \sigma^\lambda \circ \theta_{\tau ^r \land \lambda}$ for all $\lambda > 0$.
		
		Since both $(\Law(B,\tau ^\lambda))_{\lambda > 0}$ and $(\Law(\tilde{B}^\lambda,\tilde{\sigma}^\lambda))_{\lambda > 0}$ are sequences of solutions to  $\SEP(\mu,\nu)$ and $\SEP(\Law(\tilde{B}_0^\lambda),\nu)$, both families are tight by Lemma \ref{lemma:SEPcompact}. Thus, \eqref{eq:ConvPHI} holds also for all $\varphi \in C_b(\Omega \times [0, \infty))$.
		Finally, Lemma \ref{lemma:StabilityLM} shows that $\Law(\tilde{B}^\lambda,\tilde{\sigma}^\lambda)$ converges weakly $\Law(B,\tau ^{lm})$.
	\end{proof}

	\normalem

	\bibliographystyle{abbrv}
	\bibliography{lit}

\end{document}